\pdfoutput=1
\documentclass[10pt]{iopart}

\usepackage{amssymb,mathrsfs}
\usepackage{amsbsy}
\usepackage{graphicx,rotating}
\usepackage[mathscr]{euscript}
\usepackage{color}
\usepackage{rotating}
\usepackage{lscape}
\usepackage{tikz,pgfplots}
\usepackage{comment}
\usepackage{mathrsfs}
\usepackage{dsfont}
\usepackage{subfig}
\usepackage{url}

\expandafter\let\csname equation*\endcsname\relax
\expandafter\let\csname endequation*\endcsname\relax
\usepackage{amsmath}
\usepackage{mathtools}
\usepackage[pdfborder={0 0 0}, colorlinks=true, linkcolor=blue, urlcolor=blue]{hyperref}

\newcommand{\tikzw}{0.3}

\newcommand{\dom}{\mathcal{D}}
\newcommand{\dd}{\mathbf{d}}
\newcommand{\DD}{\mathbf{D}}
\newcommand{\E}{\mathcal{E}}

\newcommand{\FF}{\mathbf{F}}
\newcommand{\GG}{\boldsymbol{\Gamma}}

\newcommand{\JJ}{\mathbf{J}}
\newcommand{\Hess}{\mathcal{H}}
\newcommand{\HessGN}{\mathcal{H}_{\scriptscriptstyle\mathup{GN}}}
\renewcommand{\H}{H^1(\dom)}
\newcommand{\HH}{\mathbf{H}}

\newcommand{\Lagr}{\mathscr{L}}
\newcommand{\MM}{\mathbf{M}}

\newcommand{\nn}{\mathbf{n}}
\newcommand{\N}{\mathcal{N}} 

\renewcommand{\P}{\mathbb{P}}      
\newcommand{\R}{\mathbb{R}}       
 
\newcommand{\Us}{\mathcal{U}_s}

\newcommand{\ww}{\mathbf{w}}
\newcommand{\xx}{\mathbf{x}}
\newcommand{\yy}{\mathbf{y}}
\newcommand{\zz}{\mathbf{z}}

\newcommand{\postm}{\mu_{\rm{post}}}
\newcommand{\priorm}{\mu_0}
\newcommand{\like}{\pi_{\rm{\scriptscriptstyle{like}}} }
\newcommand{\pto}{\mathcal{F}}
\newcommand{\map}{{\scriptscriptstyle\text{MAP}}}
\newcommand{\prcov}{\mathcal{C}_0}
\newcommand{\obj}{\Phi}

\renewcommand{\SS}{\ensuremath{\mathcal{S}}}
\newcommand{\Vs}{\mathcal{V}}  

\def\addressices{Institute for Computational Engineering \& Sciences, The
  University of Texas at Austin, Austin, TX, USA}
\def\addressgeomech{Department of Geological Sciences and Department
  of Mechanical Engineering, The University of
  Texas at Austin, Austin, TX, USA}
\def\addressnyu{Courant Institute of Mathematical Sciences, New York University,
New York, NY, USA}
\def\addressncsu{Department of Mathematics, North Carolina State University, 
Raleigh, NC, USA}

\DeclareMathAlphabet{\mathup}{OT1}{\familydefault}{m}{n}

\newcommand{\Cgeneric}{\mathcal{C}_\mathup{post}}
\newcommand{\Cpost}{\mathcal{C}_\mathup{post}}
\newcommand{\Clap}{\mathcal{C}_\mathup{post}^\mathup{\scriptscriptstyle{L}}}
\newcommand{\Cgn}{\mathcal{C}^\mathup{\scriptscriptstyle{GN}}_\mathup{post}}

\begin{document}

\title[A-optimal encoding weights for nonlinear inverse problems]
{A-optimal encoding weights for nonlinear inverse problems,
with application to the Helmholtz inverse problem}

\author{Benjamin~Crestel$^1$, Alen~Alexanderian$^2$, Georg~Stadler$^3$ and
Omar~Ghattas$^{1,4}$}

\address{
$^1$\addressices\\
$^2$\addressncsu\\
$^3$\addressnyu\\
$^4$\addressgeomech
}
\ead{\url{crestel@ices.utexas.edu}, 
\url{alexanderian@ncsu.edu}, 
\url{stadler@cims.nyu.edu} and
\url{omar@ices.utexas.edu}}
\vspace{10pt}

\begin{abstract}
The computational cost of solving an inverse problem governed by PDEs, using
multiple experiments, increases linearly with the number of experiments.
A recently proposed method to decrease this cost
uses only a small number of random linear combinations of all
experiments for solving the inverse problem. This approach 
applies to inverse problems
where the PDE solution depends linearly on the right-hand side function that
models the experiment.
As this method is stochastic in essence, the quality of the obtained
reconstructions can vary, in particular when only a small number of combinations
are used.
We develop a Bayesian formulation for the
definition and computation of encoding weights that lead to a parameter
reconstruction with the least uncertainty. 
We call these weights A-optimal encoding weights.
Our framework applies to inverse problems where the governing PDE is nonlinear with
respect to the inversion parameter field. 
We formulate the problem in infinite dimensions and follow the
optimize-then-discretize approach, devoting special attention to the
discretization 
and the choice of numerical methods
in order to achieve a
computational cost that is independent of the 
parameter discretization.
We elaborate our method for a Helmholtz inverse problem, and derive
the adjoint-based expressions for the gradient of the objective function of 
the optimization problem for finding the A-optimal encoding weights.
The proposed method is potentially attractive for real-time
monitoring applications, where one can invest the effort to compute optimal
weights offline, to later solve an inverse problem repeatedly, over time,
at a fraction of the initial cost.
\end{abstract}

\vspace{2pc}
\noindent{\it Keywords}: 
source encoding,
Bayesian nonlinear inverse problem,
A-optimal experimental design,
randomized trace estimator,
Helmholtz equation.

\section{Introduction}
\label{sec:Intro}

Inverse problems are ubiquitous in science and engineering. They arise
whenever one attempts to infer 
parameters~$m$ from indirect observations~$\dd$
and from a mathematical model---the parameter-to-observable map,
$\pto(\cdot)$---for the physical phenomenon that relates $m$
and~$\dd$.
When available, it is common to use
observations obtained from different experiments to
improve the quality of the parameter estimation. Suppose
$N_s$~experiments are conduced, indexed by $i \in
\{1, \ldots, N_s \}$.  
The $i$-th experiment results in
observations~$\dd_i$ and the corresponding parameter-to-observable
map is denoted by $\pto_i(m)$.  Following a deterministic approach to this
inverse problem results in the nonlinear least-squares minimization problem
\begin{equation} \label{eq:gen} \min_m \left\{ \frac1{2N_s} \sum_{i=1}^{N_s}
\|\pto_i(m) - \dd_i \|^2 + \mathcal{R}(m) \right\}, \end{equation} where
$\mathcal{R}$ is an appropriate regularization operator to cope with
the ill-posedness that is common for many inverse problems.

Nonlinear optimization problems such as \eqref{eq:gen} can only be
solved iteratively, which requires the availability of first (and ideally,
also second) derivatives of the functional in~\eqref{eq:gen} with
respect to $m$.  For an important class of inverse problems, the
parameter-to-observable map involves the solution of a partial
differential equation (PDE). This means that the evaluation
of~$\pto_i(m)$ entails the solution~$u_i$ of a PDE, and this $u_i$
is usually restricted by an observation operator $B$ to a subset of
the domain (e.g., points), where observations are available.  In this
work, we make the assumption 
that the different experiments correspond to different right-hand
sides $f_i$ of this PDE. Moreover, this PDE must be linear with respect to the
solution~$u_i$, and both
the PDE operator as well as the observation operator $B$ must be the
same for all experiments.

When the $i$-th experiment corresponds to a forcing term~$f_i$, the
parameter-to-observable map is given by $\pto_i(m) = B u_i$, where
$\mathcal{A}(m) u_i = f_i$ with $\mathcal{A}(m)$ denoting the linear
PDE-operator that may depend nonlinearly on $m$. Note that the governing PDE can
be stationary or time-dependent. Adjoint methods allow to compute derivatives of
the objective in \eqref{eq:gen} efficiently~\cite{Troltzsch10}. For instance,
the computation of the gradient of the objective in \eqref{eq:gen} requires
solving $N_s$ forward and associated adjoint PDEs.  Similar computational costs
are associated with the application of the Hessian operator to vectors, such
that the overall computational cost of solving~\eqref{eq:gen}, which is
dominated by PDE solves with the operator $\mathcal{A}(m)$,  grows (at least)
linearly with the number of experiments $N_s$. In some important inverse
problems, $N_s$ is large (e.g., several thousand), such that these computations
are expensive or even infeasible.

There have been some recent breakthroughs to address this computational
bottleneck using the concept of random source encoding, sometimes also
referred to as simultaneous random
sources~\cite{KrebsAndersonHinkleyEtAl09,RouthLeeNeelamaniEtAl11}.  
A mathematical justification of this approach is given in the seminal
paper~\cite{HaberChungHermann12}, and is
summarized in section~\ref{sec:medparam}.  
In~\cite{LeMyersBui-Thanh16}, the authors employed a similar
idea to encode the observations in inverse problems with large
amount of data.
The main idea  of random source encoding is to replace the data
generated by each individual experiment with a small number, $N_w \ll N_s$, of linear
combinations of the data; the weights of these linear combinations, $\ww^i = [
w^i_1, \ldots, w^i_{N_s}]^T$, are called encoding weights.
Due to our linearity assumptions, this linear combination of data
corresponds to the same linear combination of experiments, i.e., we
can define encoded parameter-to-observable maps
$\pto(\ww^i;m)$, $i=1, \ldots, N_w$, as follows
\begin{equation} \label{eq:encpto}
\pto(\ww^i;m) \coloneqq \sum_{j=1}^{N_s} w^i_j \pto_j(m) = B \left(
\sum_{j=1}^{N_s} w^i_j u_j \right). \end{equation}
Observe that $\sum_{j=1}^{N_s} w^i_j u_j$ can be computed by solving the
\emph{single} PDE
\[ \mathcal{A}(m) \left( \sum_{j=1}^{N_s} w^i_j u_j \right) = \left(
\sum_{j=1}^{N_s} w^i_j f_j \right). \]  
Replacing the individual experiments with encoded experiments results
in an inverse problem with lower computational complexity.
The hope is that these linear combinations still carry
most of the information contained in the individual experiments.  
As mentioned above, the source encoding method hinges on the linearity
of the PDE describing the underlying physical phenomenon, such that
the observables depends linearly on the forcing
term. Additionally, the unicity of the observation
operator~$B$ is necessary, but this requirement can be weakened in certain situations,
e.g., if data from some experiments is missing~\cite{HaberChung14}.

The method of random source encoding, stochastic in essence, suffers from a few
limitations.
The key idea of the random source encoding approach is the
conversion of the deterministic optimization \eqref{eq:gen} into
a stochastic optimization problem. The expectation to be minimized is then
approximated using a Monte-Carlo technique (see~\cite{HaberChungHermann12} or
section~\ref{sec:medparam}). To reduce the computational cost of solving the
inverse problem, one would like to choose the number of samples used in this
Monte-Carlo approximation small. A small number of samples
translates into a large variance for the Monte-Carlo estimator of the expectation.  In
practice, this manifests itself in large differences in the reconstructions
obtained with different samples of encoding weights.
An approach to remedy that difficulty is to 
select the weights deterministically \cite{Symes10,HaberDoelHoresh14}.
In particular, in~\cite{Symes10}, the author considers to select the weights that generate the
greatest improvement from the current reconstruction, but the
results are inconclusive.
In~\cite{HaberDoelHoresh14}, the authors choose the weights that minimize the
expected medium misfit in the case of a discrete linear inverse
problem, which is related to the approach we follow in this paper.

\paragraph{Contributions} The main contributions of this article are as follows:
(1) Drawing from recent developments in optimal experimental design (OED) for high- or infinite-dimensional
inverse
problems~\cite{AlexanderianPetraStadlerEtAl14,AlexanderianPetraStadlerEtAl16,
  HaberHoreshTenorio08, HaberHoreshTenorio10},
and following a Bayesian view of inverse problems,
we develop a method for the computation of
encoding weights that lead to a parameter reconstruction with the least
uncertainty---as measured by the average of the posterior variance.
We refer to these (deterministic) weights as
\emph{A-optimal encoding weights}, a nomenclature
motivated by the use of the A-optimal experimental design criterion
from OED theory \cite{Ucinski05}.
(2) The method we propose extends the work in~\cite{HaberDoelHoresh14} by
addressing inverse problems with nonlinear parameter-to-observable maps, and
allows for infinite-dimensional parameters.  The infinite-dimensional
formulation has two main advantages: (a) the use of weak forms
facilitates the derivation of adjoint-based expressions for the gradient of the
objective function to compute the A-optimal encoding weights; (b) it allows us
to follow the optimize-then-discretize approach, which, along with devoting special attention to the
discretization of the formulation and the choice of the numerical methods
employed, helps control the computational cost independently of the
parameter discretization.
(3) We elaborate our method for the Helmholtz inverse problem and
derive the adjoint-based gradient of the optimization problem for finding the
A-optimal encoding weights.  We also analyze the computational
cost---in terms of Helmholtz PDE solves---of objective
and gradient evaluation for this optimization problem.
For this Helmholtz problem, we present an extensive numerical study and discuss the potential and
pitfalls of our approach.

\paragraph{Paper overview} The rest of this article is organized as follows.  In
section~\ref{sec:medparam}, we provide an overview of the method of
random source encoding. 
We also introduce notation that we will carry throughout the paper.
In section~\ref{sec:infdimB}, we summarize elements of Bayesian inverse
problems and introduce approximations to the posterior covariance in function
space. 
The framework for the A-optimal encoding weights is presented in
section~\ref{sec:formAweights}. 
In section~\ref{sec:wave}, we elaborate our
formulation for the Helmholtz inverse problem. We derive
adjoint-based expressions for the gradient of the A-optimal objective function,
and analyze computational cost of evaluating the objective function and its
gradient.  Numerical results are presented in
section~\ref{sec:numres}, and we provide some concluding remarks in
section~\ref{sec:conc}.

\section{Random source encoding} 
\label{sec:medparam}

In this section, we review the method of random source encoding, and introduce
notation and terminology used throughout this article.  We seek to infer a
parameter field~$m \in \Vs$ where $\Vs$ is an infinite-dimensional Hilbert space
of functions defined over the domain~$\dom \subset \R^d$ ($d=2,3$); a typical
choice is $\Vs\coloneqq L^2(\dom)$.  The parameter-to-observable map is denoted
by $\pto_i: \Vs \rightarrow \R^q$.  Let us assume that $u_i$ solves the
PDE~$\mathcal{A}(m) u_i = f_i$ and that all experiments~$i=1, \ldots, N_s$ share
a common observation operator~$B$, where $B u_i \in \R^q$. We then write each
parameter-to-observable map as $\pto_i(m) = B u_i$.  The right-hand side source~$f_i$
characterizes the $i$-th~experiment.  To apply source
encoding, we require the parameter-to-observable map to be linear with respect
to the source terms, which led us to introduce the
encoded parameter-to-observable maps~\eqref{eq:encpto}.

In~\cite{HaberChungHermann12} the authors give a mathematical justification of
the idea of random source encoding for a discrete problem and we follow their
argument, here, for an inverse problem formulated in function space.  We gather
all~$\pto_i(m)$ (resp.~$\dd_i$) in the columns of a matrix~$\FF(m)$
(resp.~$\DD^e$) and call the data misfit matrix~$\SS(m)
\coloneqq \FF(m) - \DD^e$.  Ignoring the regularization term for now, the
inverse problem can be written as, $\min_{m \in \Vs} \Big\{ \left\| \SS(m)
\right\|^2_F \Big\}$, where $\| \cdotp \|_F$ is the Frobenius
norm~\cite{TrefethenBau97}.
Note that $\| \SS(m) \|^2_F = \text{trace} ( \SS(m)^T \SS(m))$, which can
be approximated efficiently using randomized trace
estimators~\cite{AvronToledo11,Hutchinson90}.  Indeed, for random vectors~$\zz$
with mean zero and identity covariance matrix one finds that,
$\text{trace}(\SS(m)^T \SS(m)) = 
\mathbb{E}_{\zz} \big( \| \SS(m) \zz \|^2_2 \big)$.
Typical choices of distribution for~$\zz$ include the Rademacher distribution,
where samples take values~$\pm1$ with probability~$1/2$, and the standard normal
distribution~$\N(0, \mathbf{I}_{N_s})$.  Among other possible choices we mention the
discrete distribution that takes values~$\pm\sqrt{3}$ with probability~$1/6$ and
$0$ otherwise, or the uniform spherical distribution on a sphere of
radius~$\sqrt{N_s}$ that we denote~$\Us(\sqrt{N_s})$; the fact
that~$\Us(\sqrt{N_s})$
has identity covariance matrix can be shown using results
from~\cite{AndersonStephens72}, along with the observation that $\tilde{z} \sim
\Us(\sqrt{N_s})$ iff $\tilde{z} = \sqrt{N_s} (\zz / |\zz|)$ with $\zz \sim\N(0,
\mathbf{I}_{N_s})$.  We now write the data-misfit term as an expectation, i.e., $\|
\FF(m) - \DD^e \|^2_F = \mathbb{E}_{\zz} ( \| (\FF(m) - \DD^e) \zz \|^2 )$,
leading to the stochastic optimization problem
\[ \min_{m \in \Vs} \Big\{ \mathbb{E}_{\zz} \big( \| (\FF(m) - \DD^e) \zz \|^2
\big) \Big\}.  \] 
There exist two main techniques to solve these types of
problems~\cite{ShapiroDentchevaRuszczynski09}.  Using stochastic average
approximation~(SAA), one approximates the cost functional with a
Monte-Carlo-type approach before solving a deterministic optimization problem,
i.e., for fixed samples $\zz_i$ ones solves
$$\mathbb{E}_{\zz} \big( \| (\FF(m) - \DD^e) \zz \|^2 \big) \approx \frac1M
\sum_{i=1}^M \| (\FF(m) - \DD^e) \zz_i \|^2.$$ 
In an alternative approach called stochastic approximation~(SA), one re-samples
the random vector~$\zz$ at each step of the iteration.

We now specify the source-encoded equivalent of~\eqref{eq:gen}.
Given $N_w$~encoding weights~$\ww = (\ww^1, \ldots, \ww^{N_w})$, where each
$\ww^i \in \R^{N_s}$, we define the encoded data~$\dd(\ww^i) \coloneqq
\sum_{j=1}^{N_s} w_j^i \dd_j$, the encoded right-hand side~$f(\ww^i) \coloneqq
\sum_{j=1}^{N_s} w_j^i f_j$, and encoded parameter-to-observable maps
$\pto(\ww^i;m) = \sum_{j=1}^{N_s} w_j^i \pto_j(m)$.  
The parameter field~$m_c(\ww)$ reconstructed using the $N_w$~encoded sources is
then defined as
\begin{equation} \label{eq:FWILSencod} m_c(\ww) = \arg \min_{m \in \Vs} \left\{
\frac1{2N_w} \sum_{i=1}^{N_w} \left\| \pto(\ww^i;m) - \dd (\ww^i)
\right\|^2  + \mathcal{R}(m) \right\}. \end{equation} 
Due to the assumptions on~$\pto_i(m)$, the encoded map still corresponds to the
observation of a single solution to a PDE, $\pto(\ww^i;m) = B u_i$, albeit this
time $u_i$~solves the PDE $\mathcal{A}(m) u_i = f(\ww^i)$, i.e., with an encoded
right-hand side.

\section{Bayesian formulation of the inverse problem with encoded sources}
\label{sec:infdimB}

This section contains a brief presentation of the Bayesian
formulation of inverse problems with infinite-dimensional inversion
parameters; for details  we
refer the reader to~\cite{Stuart10,DashtiStuart15} for theory and to
\cite{Bui-ThanhGhattasMartinEtAl13} for the numerical approximation.
In the Bayesian framework, the unknown parameter function~$m$ is modeled as a
random field.  Starting from a prior distribution law for~$m$, we use 
observation data to obtain an improved description of the law of $m$.  This
updated distribution law of $m$ is called the posterior measure.  The prior
measure, which we denote by $\mu_0$, can be understood as a probabilistic model
for our prior beliefs about the parameter field~$m$.  The posterior measure, which we
denote by $\mu_\text{post}$, is the distribution law of $m$, conditioned on
observation data.  A key ingredient of a Bayesian inverse problem is the data
likelihood, $\like(\dd | m)$, which describes the conditional distribution of
the data given the parameter field~$m$; this is where the parameter-to-observable
map enters the Bayesian inverse problem.

Let $\dom \subset \mathbb{R}^d$ be a bounded domain with piecewise smooth
boundary and $(\Omega, \Sigma, \P)$ a probability space.  We consider an
inference parameter $m = m(x, \omega)$, with $(x,\omega) \in \mathcal{D} \times
\Omega$, such that for any $\omega \in \Omega$, $m(\cdotp, \omega) \in \Vs$
where, as before, $\Vs$ is an infinite-dimensional Hilbert space. Considering
the law of $m$ as a probability measure on $(\Vs, \mathfrak{B}(\Vs))$, the
infinite-dimensional Bayes' theorem relates the Radon-Nikodym derivative of
$\postm$ with respect to $\priorm$ with the data likelihood~$\like(\dd | m)$:
\begin{equation} \label{eq:RN} \frac{d\postm}{d\priorm} \propto \like(\dd | m) .
\end{equation} 
The use of non-Gaussian priors in infinite-dimensional Bayesian inverse problems
represents a new, interesting area of research (see for
instance~\cite{DashtiStuart15,HosseiniNigam16}). 
However, since the Bayesian inverse problem, in the formulation we introduce in
section~\ref{sec:formAweights}, only represents the
inner problem, the additional complications created by the use of
non-Gaussian priors are not justified.  We instead rely on Gaussian priors for
the Bayesian inverse problem; i.e., $\priorm = \mathcal{N}(m_0, \prcov)$ is a
Gaussian measure on $\Vs$. In that case, we require $\prcov$ to be symmetric,
positive and trace-class~\cite{Stuart10}.  A common choice for $\prcov$ (in two
and three space dimensions) is the squared inverse of a Laplacian-like
operator~$\mathcal{K}$, i.e., $\prcov = \mathcal{K}^{-2}$.  We also assume that
the noise in the data is additive, and independent and identically distributed
(over the different experiments); the distribution of each noise vector is
normal with mean zero and covariance matrix~$\GG_\text{noise}$.  That is, $\dd_i
| m \sim \mathcal{N}\big(\pto_i(m), \GG_\text{noise}\big)$, for any $i \in \{ 1,
\ldots, N_s \}$.  Consequently, each encoded observation~$\dd(\ww^i)$ will be
normally distributed with mean zero and covariance matrix ~$\GG_\text{noise,i}
\coloneqq (\sum_{j=1}^{N_s} (w_j^i)^2) \GG_\text{noise}$, i.e., $\dd(\ww^i) | m
\sim \mathcal{N}\big(\pto(\ww^i;m), \GG_\text{noise,i}\big)$, for $i \in \{ 1,
\ldots, N_w \}$.  Therefore, the likelihood function has the form 
\[ \like(\dd(\ww) | m) \propto \exp \left( - \frac1{2N_w} \sum_{i=1}^{N_w} \|
\pto(\ww^i;m) - \dd(\ww^i) \|^2_{\GG_\text{noise,i}^{-1}} \right). \]

\subsection{MAP~point}
\label{sec:MAP}
In finite dimensions, the MAP point is the parameter~$m_{\map}$ that maximizes
the posterior probability density function.  Although this definition
does not
extend directly to the infinite-dimensional case,
a MAP point can still be defined as a minimizer of a regularized data-misfit
cost functional over an appropriate Hilbert subspace of the parameter
space \cite{Stuart10}.
Let us define the Cameron-Martin space $\E = Im(\prcov^{1/2})$, endowed with the
inner-product 
\begin{equation} \label{eq:ipE}
\langle x, y \rangle_\E \coloneqq \langle \prcov^{-1/2} x, \prcov^{-1/2} y
\rangle = \langle \mathcal{K} x, \mathcal{K} y \rangle , \quad \forall x, y \in
\E. 
\end{equation}
Then the MAP~point is defined as 
\begin{equation} \label{eq:infMAPenc}
m_\map(\ww) = \arg\min_{m \in \E} \left\{ \mathcal{J}(\ww;m) \right\},
\end{equation}
where, for the inverse problems considered in the present work, 
the functional~$\mathcal{J}(\ww; \cdotp): \E \rightarrow \R$ is defined as
\begin{equation} \label{eq:Jw}
\mathcal{J}(\ww;m) \coloneqq \frac1{2N_w} \sum_{i=1}^{N_w} \left\| \pto(\ww^i;m) -
\dd(\ww^i) \right\|^2_{\GG_\text{noise,i}^{-1}}  
+ \frac12 \left\| m-m_0 \right\|^2_\E.
\end{equation}
Here, the function~$m_0 \in \E$ is the mean of the prior measure.

\subsection{Approximation to the posterior covariance}
\label{sec:approxcov}

In general, there are no closed-form expressions for moments of the posterior
measure. Thus, one usually relies on sampling-based methods to explore the
posterior. For inverse problems governed by PDEs and problems with
high-dimensional parameters (as, for instance, arising upon discretization of
an infinite-dimensional parameter field), sampling of the posterior can quickly
become infeasible since every evaluation of the likelihood requires a PDE
solve. We thus rely on approximations of the posterior, namely Gaussian
approximations about the MAP estimate. After finding the MAP point, we consider
two commonly used approximations of the posterior measure by a Gaussian measure
$\mathcal{N}(m_\map,\Cpost)$, as discussed
next~\cite{Bui-ThanhGhattasMartinEtAl13,MartinWilcoxBursteddeEtAl12}.

\paragraph{Gauss--Newton approximation} Assuming the parameter-to-observable map
$\pto(\ww^i;\cdotp)$ is Fr\'echet differentiable at the MAP point, one strategy
to approximate the posterior is to linearize around the MAP point, i.e.,
\[ \pto(\ww^i;m) \approx \pto(\ww^i;m_\map) + \JJ_{\ww^i}(m - m_\map) , \]
with $\JJ_{\ww^i}: \Vs \rightarrow \R$ the Fr\'echet derivative of the
parameter-to-observable map~$\pto(\ww^i; \cdotp)$ evaluated at the
MAP~point~\eqref{eq:infMAPenc}. Calling $(\JJ_{\ww^i})^*$ the adjoint of
$\JJ_{\ww^i}$, the covariance operator of the resulting Gaussian approximation
of the posterior is given by
\begin{equation} \label{eq:cpostgnenc} \Cgn = \left( \frac1{N_w}
\sum_{i=1}^{N_w} (\JJ_{\ww^i})^* \GG_\text{noise,i}^{-1} \JJ_{\ww^i} +
\mathcal{C}_0^{-1} \right)^{-1}.  \end{equation}
Note that the operator that appears inside the brackets in~(\ref{eq:cpostgnenc})
is the so called Gauss--Newton Hessian of the functional~(\ref{eq:Jw}) evaluated
at the MAP~point,
\[ \HessGN(m_\map) := \frac1{N_w} \sum_{i=1}^{N_w} (\JJ_{\ww^i})^*
\GG_\text{noise,i}^{-1} \JJ_{\ww^i} + \mathcal{C}_0^{-1}.  \]

\paragraph{Laplace approximation} Assuming $\mathcal{J}(\ww;\cdotp)$,
in~\eqref{eq:Jw}, is at least twice Fr\'echet differentiable at the MAP~point, a
second approach called Laplace approximation consists of using the second
derivative of $\mathcal{J}(\ww;\cdotp)$, i.e., the Hessian, at the MAP point as
an approximation to the posterior covariance 
\begin{equation} \label{eq:cpostlaplaceenc} \Clap = \left(
\mathcal{J}''(\ww;m_\map) \right)^{-1} = \Hess^{-1}(m_\map) ,  \end{equation} 
where the derivative in~$\mathcal{J}''$ is taken in terms of the parameter
field~$m$.  Note that the Laplace approximation can be related, in finite
dimensions, to a quadratic local approximation of~$\mathcal{J}(\ww; \cdotp)$
around the MAP point.

\section{A-optimal approach to source encoding}
\label{sec:formAweights}

Combining the results from section~\ref{sec:infdimB} with elements from optimal
experimental design, we propose a rigorous method to compute A-optimal
encoding weights.
In the Bayesian framework, the posterior covariance quantifies the uncertainty
in the reconstruction.  Since the posterior covariance depends on the weights
(see section~\ref{sec:postw}), we can select the weights that lead to a
reconstruction with the least uncertainty. 
In the field of optimal
experimental design, there are various design criteria that measure the
statistical quality of the reconstructed parameter field~\cite{Pukelsheim93}.  
In the present work, we rely on the A-optimal design
criterion~\cite{Pukelsheim93,AtkinsonDonev92}, which aims to minimize the trace
of the posterior covariance, 
or equivalently, to minimize the average posterior variance.
That is, we compute the weights with the smallest trace of the
posterior covariance $\obj(\ww)=\text{tr}(\Cgeneric)$, with $\Cgeneric$ given by
$\Cgn$~\eqref{eq:cpostgnenc} or $\Clap$~\eqref{eq:cpostlaplaceenc}.

An alternate view of the A-optimal design criterion is that of minimizing the
expected Bayes risk of the MAP estimator, which coincides with the trace of the
posterior covariance for a linear inverse
problem~\cite{AlexanderianPetraStadlerEtAl14,HaberHoreshTenorio08,ChalonerVerdinelli95}.
This interpretation of the A-optimal criterion can be stated as the average mean
squared error between the MAP estimator (i.e., the parameter reconstruction) and
the true parameter (e.g., see~\cite{AlexanderianPetraStadlerEtAl14}). While this
interpretation of A-optimality is restricted to linear inverse problems, it
provides another motivation for our choice of the design criterion. In our
numerical results, we explore this relation between minimizing the trace of the
posterior covariance and the mean squared distance between the MAP point and the
true parameter and observe that minimizing the trace of the posterior covariance
correlates with smaller errors for the parameter reconstruction.

\subsection{Dependence of the operators $\Cgn$ and $\Clap$ on $\ww$}
\label{sec:postw}

The dependence of the operators $\Cgn$~\eqref{eq:cpostgnenc} and
$\Clap$\eqref{eq:cpostlaplaceenc} on the weights is twofold.
First these operators depend on the encoded parameter-to-observable maps that
depend explicitly on the weights, $\pto(\ww^i;m) = \sum_{j=1}^{N_s} w_j^i
\pto_j(m)$.  Moreover, the posterior covariance operators also depend on the weights
through the MAP~point~(\ref{eq:infMAPenc}), which depends on the weights as
illustrated by~(\ref{eq:infMAPenc}) and~(\ref{eq:Jw}).

The dependence of 
the covariance operator $\Cgn$
on $\ww$ is straightforward to see.
In particular, 
using the chain-rule on the forward problem~$\mathcal{A}(m) u_i = f(\ww^i)$,
the Fr\'echet derivative of the
parameter-to-observable at the MAP~point is given by
\begin{equation} \label{eq:FF} 
\JJ_{\ww^i} = - B \mathcal{A}(m_\map(\ww))^{-1}
\frac{\partial \mathcal{A}(m) u_i}{\partial m} \Bigg|_{m = m_\map(\ww)} .
\end{equation}
Given $N_w$~encoding weights $\ww = (\ww^1, \ldots, \ww^{N_w})$ where $\ww^i \in
\R^{N_s}$, we emphasize the dependence of the posterior covariance on the
weights by writing $\Cgn = \Cgn(\ww)$.  The structure of the covariance
operator~$\Clap$ is more complicated. We detail the dependence of
$\Clap$ on $\ww$ for the application problem considered in the present paper in 
section~\ref{sec:wave}.  Note that in
the case of a linear parameter-to-observable map, both posterior
covariances~(\ref{eq:cpostgnenc}) and~(\ref{eq:cpostlaplaceenc}) are equal.

In the present formulation, $\text{tr} \big(
\mathcal{C}_\text{post}(\ww) \big)$ scales with the weights. For
instance, applying a constant multiplicative factor~$\lambda>1$ to all weights
would reduce the influence of the prior in the computation of the MAP
point~(\ref{eq:infMAPenc}) for once. It would also inflate the norm of the state
variable~$u_i$ by that factor~$\lambda$, which would then increase the size of
the derivative~(\ref{eq:FF}).  This would in turn artificially reduce the trace
of the posterior covariance~(\ref{eq:cpostgnenc}).  A solution is to restrict
the codomain of each encoding weight to a sphere of radius~$r$ in~$\R^{N_s}$.
We denote the corresponding space, for the weights~$\ww$, by~$\SS_r$, i.e., $\SS_r
\coloneqq \left\{ \ww = (\ww^1, \ldots, \ww^{N_w}) \in \R^{N_w N_s}; \, | \ww^i
| = r, \, \forall i \right\}$.  As discussed in section~\ref{sec:medparam}, the
theory of 
randomized trace estimation
dictates the use of~$r=\sqrt{N_s}$. However
this value is arbitrary and can be compensated by an equivalent re-scaling of
the regularization parameter. Therefore for simplicity we use the value~$r=1$ 
along with the notation~$\SS \coloneqq \SS_1$.
Another implication of that choice, $|\ww^i|=1$, is that the covariance matrices
for the encoded noise vectors, introduced in section~\ref{sec:infdimB}, simplify
to~$\GG_\text{noise,i} = \GG_\text{noise}$, for $i \in \{ 1, \ldots, N_w \}$.

\subsection{A-optimal encoding weights}
\label{sec:form1}

We propose to compute the A-optimal encoding weights as the solution to the
constrained minimization problem
\begin{equation} \label{eq:Wselinf}
\min_{\ww \in \SS} \obj(\ww) \coloneqq \text{tr} \big( \mathcal{C}_\text{post}(\ww) \big). 
\end{equation}
Since there are no closed-form expressions
for moments of the posterior measure, we
replace
the exact posterior covariance in~\eqref{eq:Wselinf} with one of the two
approximations introduced in section~\ref{sec:approxcov}.
The Gauss--Newton formulation of the A-optimal encoding weights,
\begin{equation} \label{eq:Wselgninf}
\obj_\mathup{GN}(\ww)  = \text{tr}(\HessGN^{-1}(\ww; m_\map(\ww)) ),
\end{equation}
is based on the posterior covariance approximation~(\ref{eq:cpostgnenc}),
and the Laplace formulation,
\begin{equation} \label{eq:Wsellaplaceinf}
\obj_\mathup{L}(\ww) = \text{tr}(\Hess^{-1}(\ww; m_\map(\ww))),
\end{equation}
is based on the posterior covariance~(\ref{eq:cpostlaplaceenc}).
Note that both formulations~\eqref{eq:Wselgninf} and~\eqref{eq:Wsellaplaceinf} require the computation of the MAP~point which
is computationally expensive for large-scale problems.  To avoid the cost
associated with the computation of the MAP~point, an additional simplification
of~\eqref{eq:Wselgninf} can be achieved  by evaluating the posterior
covariance~(\ref{eq:cpostgnenc}) at a reference parameter field~$m_0$,
which leads to the following (simplified) objective function,
\begin{equation} \label{eq:Wselgnlininf} \obj_0(\ww)  = \text{tr}(
\HessGN^{-1}(\ww; m_0)).  \end{equation}

\paragraph{A-optimal encoding weights formulation for large-scale applications}

Formulation~\eqref{eq:Wselinf} is a
nonlinear optimization problem that requires the use of iterative
methods. These methods involve repeated evaluations of the trace of the posterior
covariance.  Following discretization, the posterior covariance is a
high-dimensional  operator that is defined implicitly, i.e., through its
applications to vectors.  The exact computation of the trace of such operators,
and their derivatives with respect to encoding weights, is computationally
intractable.  For this reason, we propose an approximate formulation using a
randomized trace estimator (see~\cite{AvronToledo11,Hutchinson90} for the
theory, and~\cite{HaberDoelHoresh14,AlexanderianPetraStadlerEtAl14} for examples
of applications).
Following the formulation in~\cite{AlexanderianPetraStadlerEtAl16}, we 
introduce the Gaussian measure~$\mu_\delta = \mathcal{N}(0, \mathcal{C}_\delta)$
where $\mathcal{C}_\delta \coloneqq (I - \delta \Delta)^{-2}$.  Here $\Delta$ denotes
the Laplacian operator with homogeneous Neumann boundary conditions and $\delta
> 0$ a sufficiently small real number.  Then for any positive, self-adjoint and
trace-class operator~$\mathcal{T}$, we may use an estimator of the form, 
\[ \text{tr}(\mathcal{T}) \approx \frac1{n_{tr}} \sum_{i=1}^{n_{tr}} \langle
\mathcal{T} z_i, z_i \rangle_\mathscr{H} , \] 
where the $z_i$ are drawn from~$\mu_\delta$.  In practice, reasonable
approximations of the trace can be obtained with a relatively small $n_{tr}$.

The optimization problem for finding A-optimal encoding weights is formulated as follows 
\begin{equation*} 
\min_{\ww \in \SS} \frac1{n_{tr}} \sum_{i=1}^{n_{tr}}
\langle \mathcal{C}_\text{post}(\ww) z_i, z_i \rangle.
\end{equation*}
Specializing to the cases of
$\obj_\mathup{GN}(\ww)$~\eqref{eq:Wselgninf}
and~$\obj_{\mathup{L}}(\ww)$~\eqref{eq:Wsellaplaceinf}
results in the following 
formulations,
\begin{equation} \label{eq:Wcgninf}
\min_{\ww \in \SS} \left\{ \frac1{n_{tr}} \sum_{i=1}^{n_{tr}}
\langle \HessGN^{-1}(\ww; m_\map(\ww)) z_i, z_i \rangle \right\},
\end{equation}
\begin{equation} \label{eq:Wclaplaceinf}
\min_{\ww \in \SS} \left\{ \frac1{n_{tr}} \sum_{i=1}^{n_{tr}}
\langle \Hess^{-1}(\ww; m_\map(\ww)) z_i, z_i \rangle \right\}.
\end{equation}
Again to avoid the cost associated with the computation of the
MAP~point, one can evaluate the Gauss--Newton Hessian in~\eqref{eq:Wcgninf} at a fixed reference parameter field~$m_0$; this leads to the following (simplified) optimization problem,
\begin{equation} \label{eq:Wcgnlininf} 
\min_{\ww \in \SS} \left\{ \frac1{n_{tr}} \sum_{i=1}^{n_{tr}}
\langle \HessGN^{-1}(\ww; m_0) z_i, z_i \rangle \right\}.
\end{equation}
The formulation~\eqref{eq:Wcgnlininf} can be seen as an extension of the
formulation proposed in~\cite{HaberDoelHoresh14} to a fully nonlinear inverse
problem formulated at the infinite-dimensional level.

\section{Application to the Helmholtz inverse problem}
\label{sec:wave}

In this section, we elaborate 
the A-optimal encoding weights formulation introduced in
section~\ref{sec:formAweights} for the 
Helmholtz inverse problem.
Recall that high resolution
reconstructions in this application require a large number of experiments and
that the computational cost of the inversion scales linearly with the
number of experiments (see section~\ref{sec:Intro}).  Source encoding can
provide a trade-off between high-quality reconstruction and computational cost.

We begin by describing the inverse problem used in our study (section
\ref{sec:inverse_problem}).  
Then the optimization problem to compute the A-optimal encoding weights,
including the adjoint-based expressions for the gradient of this
objective function, is detailed in section~\ref{sec:optimal_weights}.  

\subsection{The inverse problem: medium parameter reconstruction} 
\label{sec:inverse_problem}

For simplicity of
the presentation, we derive the formulation using a single frequency but 
extensions to the case of multiple frequencies are straightforward.
We use homogeneous Neumann boundary conditions.  The frequency-domain
Helmholtz equation is given, for~$i = 1,\ldots, N_w$, by
\begin{equation} \label{eq:strongH} 
\begin{aligned} 
   -\Delta u_i - \kappa^2 m u_i & = f(\ww^i)
  & \text{in } \dom, \\ 
\nabla u_i \cdotp \nn & = 0
 & \text{on } \partial \dom.  \end{aligned} 
\end{equation} 
Solutions~$u_i$ \eqref{eq:strongH} are considered in $\H$, i.e.,
the Sobolev space of functions in $L^2(\dom)$ with square integrable weak
derivatives.
The original source terms are in the dual space of~$H^1_0(\dom)$, i.e.,~$f_j \in
H^{-1}(\dom)$.
The (medium) parameter field~$m \in L^\infty(\dom)$ 
corresponds to the square of the slowness (or the squared
inverse local wave speed) and the constant $\kappa$ is the frequency of the wave
(in rad/s).

\subsubsection{MAP point} 
\label{sec:MAPpt}

The MAP~point is the solution to a deterministic inverse problem
(see section~\ref{sec:MAP}) with the norms in the data-misfit and regularization terms
weighted by the noise and prior covariance operators respectively.
In particular, with a Gaussian prior~$\priorm =
\mathcal{N}(m_0, \mathcal{C}_0)$ and the norm corresponding to the inner
product~\eqref{eq:ipE}, we have 
\begin{equation} \label{eq:invpb} m_\map(\ww) = \arg \min_{m \in \E} \left\{
\frac1{2N_w}  \sum_{i=1}^{N_w} \left\| B u_i - \dd(\ww^{i})
\right\|^2_{\GG^{-1}_\text{noise}} + \frac12 \left\| m-m_0 \right\|^2_\E
\right\} , \end{equation} 
where $u_i$ solves~(\ref{eq:strongH}).  

To properly define the source terms $f_i$, appearing in the right hand-side of
the forward problem, and the observation operator $B$, we define the
mollifier~$\varphi_\varepsilon(x; y)$ as follows: 
\begin{equation} \label{eq:mollifier}
\varphi_\varepsilon(x;y) = \frac1{\alpha_\varepsilon}
e^{-\frac1{\varepsilon^2 - |x-y|^2}}
\mathds{1}_{\mathcal{B}(y,\varepsilon)}(x) , \end{equation}
where $\alpha_\varepsilon = 2\pi K \varepsilon^2 e^{-1/\varepsilon^2}$, $K =
\int_0^1 r e^{-1/(1-r^2)} dr$, $\mathds{1}_{\mathcal{B}(y,\varepsilon)}$ is the
indicator function for the ball of radius~$\varepsilon$ centered at~$y$,
and $0 < \varepsilon \ll 1$. This function is smooth and integrates to one.  
We choose each source terms~$f_i$ to be a mollifier centered at one
of the $N_s$~source locations that we denote~$x_i^s$ for $i=1, \ldots, N_s$,
i.e., $f_i(x) = \varphi_\varepsilon(x; x_i^s)$.  The observation operator~$B: \H
\rightarrow \R^q$ is the evaluation, at each of the receiver locations which we
denote~$x_j^r$ for $j=1, \ldots, q$, of a convolution between the solution to
the forward problem~$u_i$ and a mollifier~$\varphi_{\varepsilon'}(x;0)$, i.e.,
$(B u_i)_j = (u_i * \varphi_{\varepsilon'}(\cdotp;0))(x_j^r)$.  These choices of
the source terms and observation operator guarantee that the forward, adjoint,
incremental forward and incremental adjoint solutions belong to~$\H$.

\subsubsection{Gradient and Hessian of the inverse problem} 

Availability of derivatives of the function in brackets on the right hand side
of \eqref{eq:invpb} is required for the computation of $m_\map$.  The second
derivative, i.e., the Hessian operator, also enters the A-optimal formulation
laid down in section~\ref{sec:formAweights}.
We derive both gradient and
Hessian following the formal Lagrangian approach~\cite{Troltzsch10,BorziSchulz12}. 
The first-order
necessary optimality condition for the MAP~point is 
a coupled system of PDEs: Find $(m_\map,\{u_i\}_i,\{p_i\}_i) \in \E
\times \H^{N_w} \times \H^{N_w}$ such that for all variations $(\tilde{m},
\{\tilde{u}_i\}_i, \{ \tilde{p}_i \}_i) \in \E \times \H^{N_w} \times \H^{N_w}$
\begin{equation} \label{eq:fooc} \begin{aligned} 
\langle \nabla u_i , \nabla \tilde{p}_i \rangle - \kappa^2 \langle m_\map(\ww) u_i ,
\tilde{p}_i \rangle - \langle f(\ww^i), \tilde{p}_i \rangle  &
= 0, \, \forall i \\
\langle \nabla \tilde{u}_i, \nabla p_i \rangle - \kappa^2 \langle \tilde{u}_i,
m_\map(\ww) p_i \rangle +\langle B\tilde{u}_i, B  u_i - \dd(\ww^i)
\rangle_{\GG^{-1}_\text{noise}} & =0, \, \forall i \\ 
\langle m_\map(\ww) -m_0, \tilde{m} \rangle_\E -\frac1{N_w} \sum_{i=1}^{N_w} \kappa^2
\langle u_i p_i, \tilde{m} \rangle & =0.  
\end{aligned} \end{equation}
For the Hessian, we describe the solution to the equation~$y =
\Hess^{-1}(m_\map) z$.  This leads to the coupled system of PDEs: Find
$(y,\{v_i\}_i,\{q_i\}_i) \in \E \times \H^{N_w} \times \H^{N_w}$ such that for
all $(\tilde{m}, \{\tilde{u}_i\}_i, \{ \tilde{p}_i \}_i) \in \E \times \H^{N_w}
\times \H^{N_w}$ the following equations are satisfied:
\begin{equation} \label{eq:Hz} \begin{aligned} 
\langle \nabla v_i, \nabla \tilde{p}_i \rangle
- \kappa^2 \langle m_\map(\ww) v_i, \tilde{p}_i \rangle
- \kappa^2 \langle u_i y, \tilde{p}_i \rangle & = 0 , \, \forall i\\ 
\langle \nabla
  \tilde{u}_i, \nabla q_i \rangle
- \kappa^2 \langle\tilde{u}_i, m_\map(\ww) q_i \rangle
- \kappa^2 \langle \tilde{u}_i, p_i y \rangle + \langle B \tilde{u}_i,
   B v_i \rangle_{\GG^{-1}_\text{noise}} & = 0 , \, \forall i \\ 
\langle y, \tilde{m} \rangle_\E
- \frac1{N_w} \sum_{i=1}^{N_w} \kappa^2 \Big[ \langle v_i p_i, \tilde{m} \rangle
  + \langle  u_i q_i, \tilde{m} \rangle \Big] & = \langle z, \tilde{m}
\rangle.  \end{aligned} \end{equation}

\subsection{The optimization problem for A-optimal encoding weights} 
\label{sec:optimal_weights}

Here we formulate the optimization problem for computing A-optimal source
encoding weights for the frequency-domain seismic inverse problem
\eqref{eq:strongH}. We restrict ourselves to the case of the Laplace
formulation~(\ref{eq:Wclaplaceinf}) as the other two functionals,
(\ref{eq:Wcgninf}) and~(\ref{eq:Wcgnlininf}), can be treated as special cases
of the Laplace formulation.

In its original format, the optimization problem for A-optimal encoding
weights~\eqref{eq:Wclaplaceinf} is a bi-level optimization, as the MAP~point is
itself the solution to a minimization problem~\eqref{eq:infMAPenc}. 
However this is not a practical formulation to compute derivatives. 
We therefore reformulate~\eqref{eq:Wclaplaceinf} as a PDE-constrained
optimization problem in which the MAP point is defined as a solution of the
first-order optimality condition~\eqref{eq:fooc}.  The other PDE constraint is
the solution to the Hessian system~\eqref{eq:Hz} along the random directions of
the trace estimator, i.e., we define the objective functional for the
computation of the A-optimal encoding weights by
\[ \frac1{n_{tr}} \sum_{k=1}^{n_{tr}} \langle y_k, z_k \rangle, \]
where $z_k$ is a random direction for the trace estimator and $y_k =
\Hess^{-1}(m_\map) z_k$ according to~\eqref{eq:Hz}.  We can then enforce these
PDE constraints with Lagrange multipliers and compute derivatives of the
optimization problem~\eqref{eq:Wclaplaceinf} using the formal
Lagrangian approach.
We account for the constraint on the weights through a penalty term,
\[ \frac{\lambda}{2 N_w} \sum_{j=1}^{N_w} \left( \| \ww^j \|^2 - 1 \right)^2 , \] 
with $\lambda \in \R$.  Although a penalty term is not the only option, we found
this relaxation of the constraint to be efficient and easy to implement.

We now present the complete formulation for~\eqref{eq:Wclaplaceinf}.  The
A-optimal encoding weights are solutions to the minimization problem
\begin{equation} \label{eq:wkform} \min_\ww \Bigg\{ \frac1{n_{tr}}
\sum_{k=1}^{n_{tr}} \langle y_k, z_k \rangle + \frac{\lambda}{2 N_w}
\sum_{j=1}^{N_w} \left( \| \ww^j \|^2 - 1 \right)^2 \Bigg\} ,\end{equation}
where for every $k=1, \ldots, n_{tr}$, 
$(y_k,\{v_{i,k}\}_i,\{q_{i,k}\}_i) \in \E \times \H^{N_w} \times
\H^{N_w}$ solves the system
\begin{equation} \label{eq:costinvH} \begin{aligned} 
\langle \nabla v_{i,k}, \nabla \tilde{p}_{i,k} \rangle - \kappa^2 \langle m_\map(\ww)
v_{i,k}, \tilde{p}_{i,k} \rangle - \kappa^2 \langle u_{i} y_k, \tilde{p}_{i,k}
\rangle & = 0 , \, \forall i \\ 
\langle \nabla \tilde{u}_{i,k}, \nabla q_{i,k} \rangle - \kappa^2
\langle\tilde{u}_{i,k}, m_\map(\ww) q_{i,k} \rangle - \kappa^2 \langle
\tilde{u}_{i,k}, p_{i} y_k \rangle \hspace{.5in} & \\
+ \langle B \tilde{u}_{i,k},
 B v_{i,k} \rangle_{\GG^{-1}_\text{noise}} & = 0 , \, \forall i \\ 
\langle y_k, \tilde{m} \rangle_\E - \frac1{N_w} \sum_{i=1}^{N_w} \kappa^2 \Big[
\langle v_{i,k} p_{i}, \tilde{m} \rangle   + \langle  u_{i} q_{i,k}, \tilde{m}
\rangle \Big] & = \langle z_k, \tilde{m} \rangle,  
\end{aligned} \end{equation}
for all $(\tilde{m}, \{\tilde{u}_{i,k}\}_i, \{ \tilde{p}_{i,k} \}_i) \in \E
\times \H^{N_w} \times \H^{N_w}$ and where $(m_\map,\{u_i\}_i,\{p_i\}_i) \in \E
\times \H^{N_w} \times \H^{N_w}$ solves the first-order optimality system for
the Helmholtz inverse problem
\begin{equation*} \begin{aligned} 
\langle \nabla u_i  , \nabla \tilde{p}_i \rangle - \kappa^2 \langle m_\map(\ww)
u_i  , \tilde{p}_i  \rangle - \langle f(\ww^i) , \tilde{p}_i  \rangle & = 0, \,
\forall i \\ 
\langle \nabla \tilde{u}_i , \nabla p_i  \rangle - \kappa^2 \langle \tilde{u}_i
, m_\map(\ww) p_i  \rangle +\langle B\tilde{u}_i , B  u_i  - \dd(\ww^i)
\rangle_{\GG^{-1}_\text{noise}} & =0, \, \forall i \\ 
\langle m_\map(\ww)-m_0, \tilde{m} \rangle_\E -\frac1{N_w} \sum_{i=1}^{N_w}
\kappa^2 \langle u_i  p_i , \tilde{m} \rangle & =0, 
\end{aligned} \end{equation*}
for all $(\tilde{m}, \{\tilde{u}_i\}_i, \{ \tilde{p}_i \}_i) \in \E \times
\H^{N_w} \times \H^{N_w}$.

\subsubsection{Gradient of the A-optimal weight problem} \label{sec:gradout} 

We derive the gradient of the objective function defined in~\eqref{eq:wkform},
with respect to $\ww$, using a formal Lagrangian approach.  
We refer the reader to~\ref{sec:outer} for this derivation.  
Since we enforce the PDE constraints weakly using Lagrange multipliers,
we introduce adjoint
variables that are indicated with a star superscript, e.g., $m^*$ is the adjoint
variable for~$m$.  
The gradient is given by $\left[ \delta_{\ww^1} \obj_L(\ww), \delta_{\ww^2}
\obj_L(\ww), \ldots, \delta_{\ww^{N_w}} \obj_L(\ww) \right]^T$, where for
any~$i=1,\ldots, N_w$,
\begin{equation*} 
\delta_{\ww^i} \obj_L(\ww) =   - \frac1{N_w} 
\begin{bmatrix} 
\langle f_1, u_i^* \rangle & + & \langle B p_i^*, \dd_1 \rangle_{\GG^{-1}_\text{noise}}\\ 
\langle f_2, u_i^* \rangle & + & \langle B p_i^*, \dd_2 \rangle_{\GG^{-1}_\text{noise}} \\ 
& \vdots & \\ 
\langle f_{N_s}, u_i^* \rangle & + & \langle B p_i^*, \dd_{N_s} \rangle_{\GG^{-1}_\text{noise}}
\end{bmatrix}. 
\end{equation*}

The variables~$u_i^*$ and $p_i^*$ are computed by solving the following
Hessian-like system (compare with~\eqref{eq:Hz}): Find
$(m^*,\{u_i^*\}_i,\{p_i^*\}_i) \in \E \times \H^{N_w} \times \H^{N_w}$ such that
for all $(\tilde{m}, \{\tilde{u}_i\}_i, \{ \tilde{p}_i \}_i) \in \E \times
\H^{N_w} \times \H^{N_w}$ the following equations are satisfied: 
\begin{equation} \label{eq:mstar2} \begin{aligned}
\langle \nabla p_i^*, \nabla \tilde{p}_i \rangle
- \kappa^2 \langle m p_i^*, \tilde{p}_i \rangle 
- \kappa^2 \langle u_i m^*, \tilde{p}_i \rangle & = - \frac2{n_{tr}}
  \sum_{k=1}^{n_{tr}} \kappa^2 \langle v_{i,k} y_k , \tilde{p}_i \rangle, \\ 
\langle \nabla u_i^*, \nabla \tilde{u}_i  \rangle 
- \kappa^2 \langle m u_i^*, \tilde{u}_i \rangle 
- \kappa^2 \langle p_i m^*, \tilde{u}_i \rangle \hspace{.2in} & \\ 
 + \langle B p_i^*, B \tilde{u}_i \rangle_{\GG^{-1}_\text{noise}}  & = -
\frac2{n_{tr}} \sum_{k=1}^{n_{tr}} \kappa^2 \langle y_k q_{i,k}, \tilde{u}_i \rangle,
\\
\langle m^*, \tilde{m} \rangle_\E
- \frac1{N_w} \sum_{i=1}^{N_w} \kappa^2 \left[ \langle u_i u_i^*, \tilde{m} \rangle +
  \langle p_i^*p_i,\tilde{m} \rangle \right] & =
- \frac2{n_{tr}N_w} \sum_{k=1}^{n_{tr}} \sum_{i=1}^{N_w} \kappa^2 \langle v_{i,k}
  q_{i,k}, \tilde{m} \rangle .  \end{aligned} \end{equation}
The variables~$\{v_{i,k}\}$ (resp.~$\{q_{i,k}\}$) are the incremental state
(resp.~adjoint) variables which occur in the application of the inverse Hessian in the
direction of the $k$-th~trace estimator direction~$z_k$.  

\subsubsection{Discretization} 
\label{sec:discrete}

The numerical solution of~(\ref{eq:wkform}) is done via the
Optimize-then-Discretize~(OTD) approach, where the discretization is based on
continuous Galerkin finite element with Lagrange nodal basis functions.  
Extra care is needed for the discretization of the
covariance operator to ensure that its discrete representation
faithfully represents the properties of the target infinite-dimensional object.
We do not provide full details of the discretization and refer the reader
to~\cite{AlexanderianPetraStadlerEtAl14,Bui-ThanhGhattasMartinEtAl13}.  However,
we show how to select the discrete random directions~$z_k$ in the trace
estimator.  Let us call~$V_h$ the finite-dimensional approximation to the space
$\H$ used for the finite-element representations of all state, adjoint,
corresponding incremental variables and their respective adjoints.  And let
$V_h^m$ be the finite-dimensional space for the medium parameter~$m$.  Let us
call $\{ \psi_i\}_{i=1}^t$ (resp.~$\{ \phi_i \}_{i=1}^l$) a basis for~$V_h$
(resp.~$V_h^m$).  Let us introduce the vector notations $\yy_k = ( y_k^1,
\ldots, y_k^l)^T$ (resp.~$\zz_k = (z_k^1, \ldots, z_k^l )^T$) for the
finite element representations of $y_k$ (resp.~$z_k$) in~$V_h^m$.  The
finite-dimensional approximation to the trace estimation is then
\[ \frac1{n_{tr}} \sum_{k=1}^{n_{tr}} \langle y^h_k, z^h_k \rangle_{L^2} =
\frac1{n_{tr}} \sum_{k=1}^{n_{tr}} \sum_{i,j=1}^l y_k^i z_k^j \langle \phi_i,
\phi_j \rangle_{L^2}  = \frac1{n_{tr}} \sum_{k=1}^{n_{tr}} \langle \yy_k, \zz_k
\rangle_\MM , \] 
with $\MM_{ij} = \langle \phi_i, \phi_j \rangle_{L^2}$ the mass matrix
in~$V_h^m$.  From the definition of~$y_k$, we see that each $y_k^h$ solves the
system $\langle \Hess y_k^h, \phi_i \rangle_{L^2} = \langle z_k^h, \phi_i
\rangle_{L^2}$, for $i=1, \ldots,l$.  Substituting the representation of~$y_k^h$
and~$z_k^h$ in the basis of~$V_h^m$, we obtain the matrix system $\HH \yy_k =
\MM \zz_k$, where $\HH$ is the standard Hessian matrix obtained from
finite-element discretization of system~(\ref{eq:Hz}), i.e., $\HH_{ij} = \langle
\Hess \phi_j, \phi_i \rangle_{L^2}$.  The finite-dimensional approximation to
the trace estimation becomes
\[ \frac1{n_{tr}} \sum_{k=1}^{n_{tr}} \langle y^h_k, z^h_k \rangle_{L^2}  =
\frac1{n_{tr}} \sum_{k=1}^{n_{tr}} \langle \HH^{-1} \MM \zz_k, \zz_k \rangle_\MM
= \frac1{n_{tr}} \sum_{k=1}^{n_{tr}} \langle \HH_\MM^{-1} \zz_k, \zz_k
\rangle_\MM, \]
where we defined $\HH^{-1}_\MM \coloneqq\HH^{-1} \MM$.  The matrix
$\HH^{-1}_\MM$ is $\MM$-symmetric~\cite{Bui-ThanhGhattasMartinEtAl13}, i.e.,
self-adjoint with respect to the $\MM$ inner-product. Then it was proved
in~\cite{AlexanderianPetraStadlerEtAl14} that $\frac1{n_{tr}}
\sum_{k=1}^{n_{tr}} \langle \HH_\MM^{-1} \zz_k, \zz_k \rangle_\MM$ is indeed a
trace estimator provided $\zz_k \sim \mathcal{N}(0, \MM^{-1})$.  In practice,
vectors $\zz_k$ are sampled by taking draws $\xx_k$ 
from multivariate standard normal distribution, $\xx_k \sim \mathcal{N}(0, \mathbf{I})$,
and using $\zz_k = \MM^{-1/2} \xx_k$ 

\subsubsection{Computational cost}\label{sec:computational_cost}
\newcommand{\nCG}{n_\text{cg}}
\newcommand{\nNw}{n_\text{newt}}

Problem~\eqref{eq:wkform} is highly nonlinear and requires iterative methods to
be solved. The gradient, derived in section~\ref{sec:gradout}, allows us to
use quasi-Newton methods~\cite{NocedalWright06}.  
In table~\ref{tab:cost}, we 
report the dominating terms of the computational
cost of evaluating the objective function and its gradient in all three
cases~\eqref{eq:Wcgninf}-\eqref{eq:Wcgnlininf}.
Additionally, it is possible to
reduce the cost of formulation~\eqref{eq:Wcgninf} by computing a low-rank
approximation of the Hessian operator~\cite{FlathWilcoxAkcelikEtAl11}.  One must
keep in mind, however, that the incremental state variables~$\{ v_{i,k} \}$ and
incremental adjoint variables~$\{ q_{i,k} \}$ corresponding to each random
directions~$\{ z_k \}$ are required to compute the gradient.
For this reason, a low-rank approximation of the Hessian will only lower the
computational cost when $n_{tr} > \nCG \nNw$.

\begin{table}[h]
\centering
\caption{Computational cost 
for objective function and gradient evaluation of the optimization problem
for finding A-optimal encoding weights. We report the computational
cost, in terms of the number of forward PDE solves, for 
$\obj_{\mathup{GN}}(\ww)$, $\obj_{\mathup{L}}(\ww)$, and $\obj_0(\ww)$
defined in~\eqref{eq:Wcgninf}--\eqref{eq:Wcgnlininf} respectively. 
Notations: $\nCG=$~number of Conjugate-Gradient iterations
to compute the search direction in Newton's method; $\nNw=$~number of Newton
steps to compute the MAP~point.}
\begin{tabular}{ll|c|c|c}
\multicolumn{2}{c|}{ } 
& $\!\!\!\! \obj_0(\ww)\!$
$\!\!\!\!\!$
& $\!\!\obj_\mathup{GN}(\ww)\!$ $\!~${\small and}$\!$
$ \obj_\mathup{L}(\ww)$ $\!\!\!\!\!$ 
& $\!\! \obj_\mathup{GN}(\ww)\!$ $\!\!$ \\
 & & & (no low-rank) & (with low-rank) \\ \hline\hline
\multicolumn{2}{l|}{objective evaluation} & & & \\
 & MAP~point & $2N_w$ & $\!\!2N_w\nCG\nNw\!\!$ & $\!\!2N_w\nCG\nNw\!\!$ \\
 & tr$(\Hess^{-1})$ & $\!\! 2N_w\nCG n_{tr} \!\!$ & $2N_w\nCG n_{tr}$ & $2N_w\nCG$ \\
\hline
\multicolumn{2}{l|}{gradient evaluation} & & & \\
 & $v_{ik}$, $q_{ik}$ & --  & -- & $2 N_w n_{tr}$ \\
 & $m^*$ & -- & $2 N_w \nCG$ & -- \\
 & $u_i^*$, $p_i^*$ & $N_w$  & -- & $2N_w$ \\
\hline
\hline
\multicolumn{2}{l|}{total \!\!\!\!\!} & $\!\!\!\! 2N_w \nCG n_{tr}
\!\!\!\!$ &
\multicolumn{1}{|c|}{$\!\! 2N_w \nCG (\nNw \!+\! n_{tr})\!\!$} &
\multicolumn{1}{|c}{$\! 2N_w (\nCG\nNw \!+\! n_{tr}) \!\!\!$} \\  
\end{tabular}
\label{tab:cost}
\end{table}

Following the OTD approach, the optimization problem~\eqref{eq:wkform} is
formulated in function space, before being solved with algorithms that are
discretization-independent.  
This results in the overall computational cost
being independent of the discretization of the parameter space, or in other
words, each of the quantities $\nNw$, $\nCG$ and $n_{tr}$ in
table~\ref{tab:cost} remain constant when the mesh gets refined. 
We spend the
rest of this section discussing the choice of such discretization-invariant
algorithms.
First, we use Newton's method, with Armijo line search, to compute the
MAP~point; the number of Newton steps needed to converge,~$\nNw$, is typically independent
of the size of the parameter space~\cite{Deuflhard04}.
Moreover, the Hessian system~\eqref{eq:Hz} needed to compute
the MAP~point, to evaluate the objective function~\eqref{eq:wkform}, and to compute
the adjoint variable~$m^*$~\eqref{eq:mstar2}, is solved using the preconditioned
Conjugate Gradient method~\cite{NocedalWright06}.
The Conjugate Gradient solver is
preconditioned by the prior covariance operator; the number of iterations~$\nCG$
needed to solve the Hessian system then depends on the spectral properties of
the prior-preconditioned data-misfit part of the Hessian operator (i.e., the
Hessian in function space) and is therefore independent of the discretization. 
The trace estimator displays a similar type of behaviour.  The number of trace
estimator vectors~$n_{tr}$ one should use depends on the spectral properties of the
underlying infinite-dimensional operator.  The choice of a discrete
inner-product weighted by the mass matrix (see section~\ref{sec:discrete}) guarantees that our discrete
operator will be a valid approximation of the infinite-dimensional operator and
will conserve its spectral properties. The actual evaluation of the trace is
performed through the repeated solution of the Hessian
system~\eqref{eq:costinvH}, which was shown above to be
discretization-independent.

\section{Numerical results}
\label{sec:numres}

In this section, we present  numerical results for the Helmholtz inverse problem
in two (spatial) dimensions.  We start with a low-dimensional example ($N_w=1$
for $N_s=2$), which allows us to visualize the objective functions defined in
section~\ref{sec:form1} over the entire weight space.  This facilitates a
qualitative comparison of the different approximations introduced, the
Gauss--Newton~\eqref{eq:Wselgninf} and Laplace objective
functions~\eqref{eq:Wsellaplaceinf}, along with the linearized
formulation~\eqref{eq:Wselgnlininf}.  We then present an example with a
higher-dimensional weight space~($N_s=10$) in which we study the distribution of
the A-optimal encoding weights and random weights sampled from the uniform
spherical distribution and how the number of encoded weight vectors influence
these results.

The setting for this section is a square domain with 20 receivers located at the
top of the domain, and sources positioned on the bottom and left edges of the
domain.  The source term is a mollifier~\eqref{eq:mollifier} with
$\varepsilon=10^{-6}$.  This choice of source terms was numerically found to be
reasonably well approximated, at the discrete level, by a point source; we
utilize that approximation in this section.
We use a wave frequency of~$\kappa = 2\pi$ in equation~\eqref{eq:strongH}.  All
partial differential equations are discretized by continuous Galerkin finite
elements (linear elements for the parameters and quadratic elements for the
state and adjoint variables). This results in a (medium) parameter space of
182~degrees of freedom. We work with synthetic data that are polluted by a 2\%
additive Gaussian noise.

\subsection{One-dimensional weight space}
\label{sec:1dspace}

\begin{figure}
\centering
\begin{tabular}{ccc}
\begin{tikzpicture}
\node[anchor=south west, inner sep=0] (image1) at (0,0)
{\includegraphics[height=0.25\textwidth, trim=50 20 90 10,
clip=true]{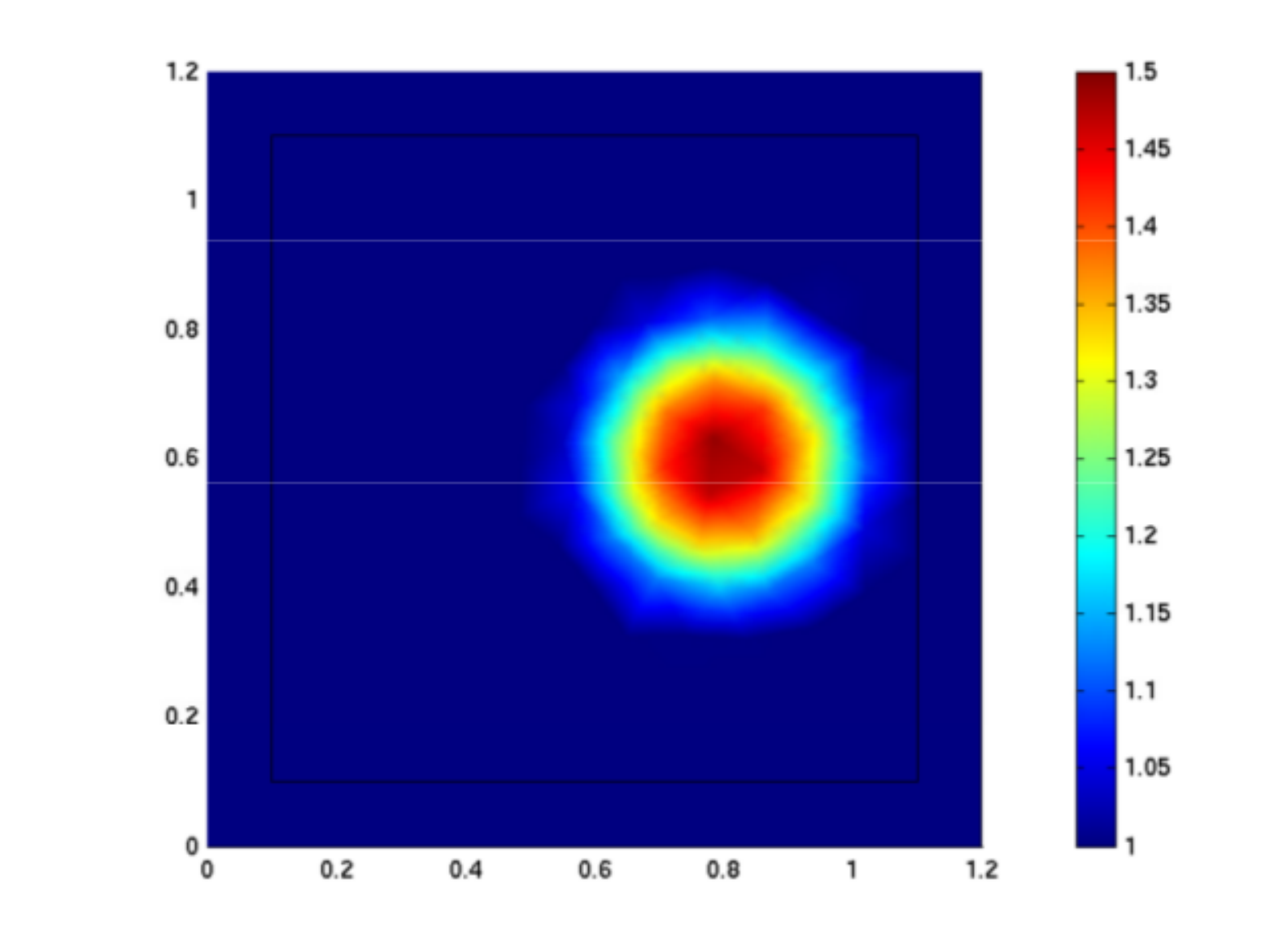}};
\begin{scope}[x={(image1.south east)},y={(image1.north west)}]
\filldraw[fill=green,draw=black] (0.05,0.17) rectangle (0.11,0.23);
\filldraw[fill=green,draw=black] (0.34,0.03) rectangle (0.4, 0.09);
\foreach \x in {0,1,...,19} {
\filldraw[fill=yellow,draw=black] (0.15+0.03895*\x,.95) circle (0.02);}
\end{scope}\label{fig:domsrc1}
\end{tikzpicture}
&
\includegraphics[height=0.25\textwidth, trim=350 20 30 10,
clip=true]{targetmedium-1}
&
\begin{tikzpicture}
\node[anchor=south west, inner sep=0] (image1) at (0,0)
{\includegraphics[height=0.25\textwidth, trim=50 20 90 10,
clip=true]{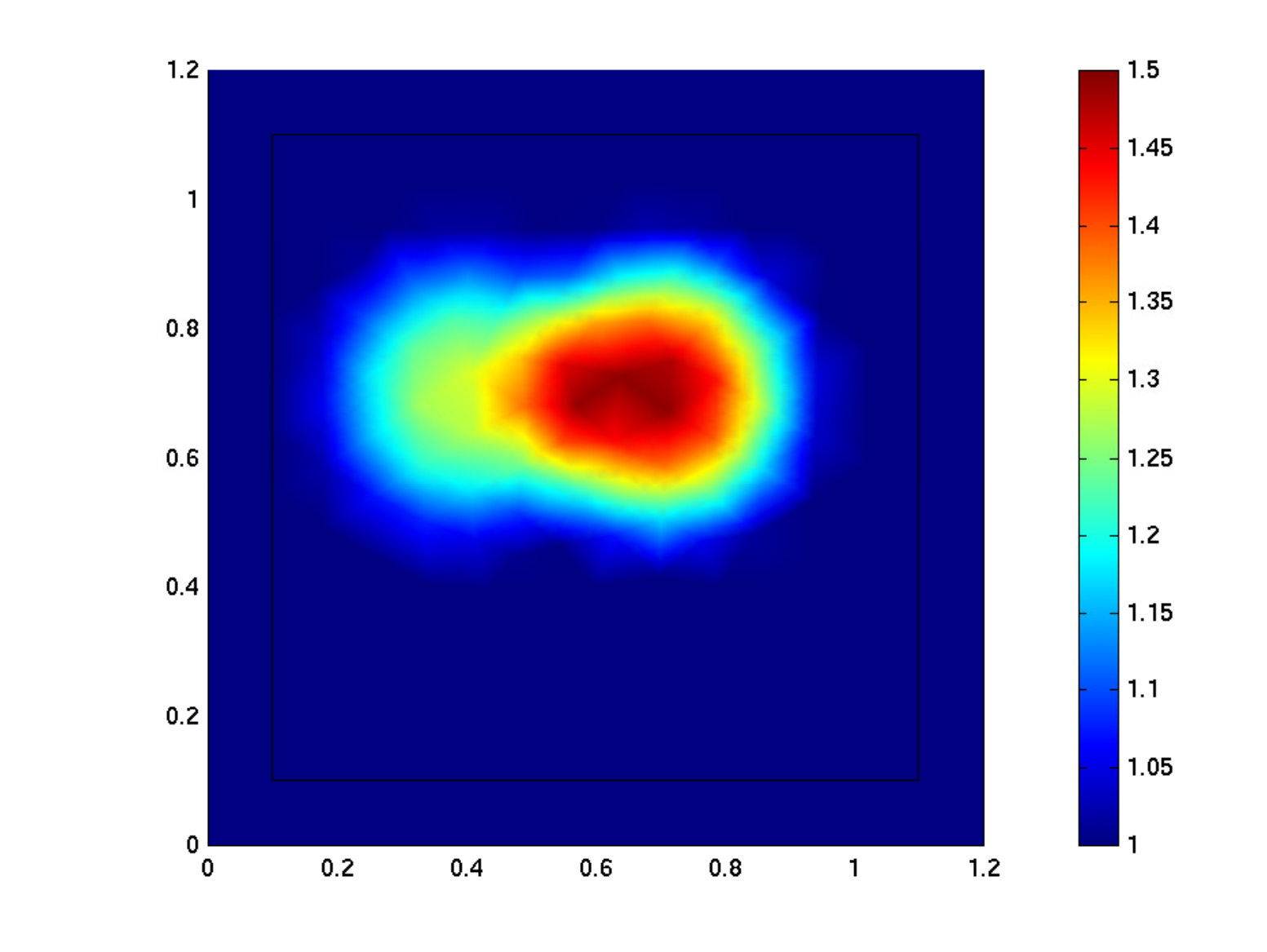}};
\begin{scope}[x={(image1.south east)},y={(image1.north west)}]
\filldraw[fill=green,draw=black] (0.05,0.17) rectangle (0.11,0.23);
\filldraw[fill=green,draw=black] (0.34,0.03) rectangle (0.4, 0.09);
\foreach \x in {0,1,...,19} {
\filldraw[fill=yellow,draw=black] (0.15+0.03895*\x,.95) circle (0.02);}
\end{scope}\label{fig:domsrc2}
\end{tikzpicture}
\\
(a) medium parameter 1 & & (b) medium parameter 2
\end{tabular}
\caption{Target medium parameters, along with the locations of the sources
(green squares) and receivers (yellow circles). 
}
\label{fig:dom2src}
\end{figure}

In this section, we study a one-dimensional source encoding problem
corresponding to a single linear combination of two sources ($N_s=2$ and
$N_w=1$).  Although this setting represents an unrealistic situation (low number
of sources, and high ratio of number of encoded sources over total number of
sources), it is informative for the following reasons: (1) It provides numerical
evidence of the strong and highly nonlinear dependence of the objective
functions~\eqref{eq:Wselgninf}--\eqref{eq:Wselgnlininf} on the encoding weights.
(2) It demonstrates the presence of multiple local minima in the minimization
problem~\eqref{eq:Wselinf}. (3) It highlights the difference between the
Gauss--Newton and Laplace formulations.  The sources are located on the bottom
and left edges of the domain, and we study two different medium parameters, each
made of a constant background and a smooth compactly supported perturbation (see
figure~\ref{fig:dom2src}).

We next define the noise covariance and the prior covariance operators used in
these numerical applications.  Let us introduce the non-singular, positive
definite, elliptic operator~$\mathcal{Y}= -\gamma \Delta + \beta I$, with
$\gamma, \beta$ positive constants, $I$ the identity operator and $\Delta$ the
Laplacian operator with homogeneous Neumann boundary conditions. Then we define
the prior covariance operator as $\prcov^{-1} = \mathcal{Y} + \eta
\mathcal{Y}^2$ with $\eta > 0$.  One can verify that this choice of prior
covariance operator is symmetric, positive definite and trace-class as long as
$\gamma, \eta, \beta > 0$.  The noise covariance operator for the observations
is chosen to be a multiple of the identity matrix, i.e., $\GG_\text{noise} =
\sigma^2 \mathbf{I}$---in our examples we choose $\sigma=1$.  The
parameters~$\gamma$, $\beta$, and~$\eta$ are chosen as $\gamma = 10^{-3}$,
$\beta = 10^{-4}$ and $\eta = 10^{-2}$, and we have verified that this choice
approximately satisfies the discrepancy principle.  In the (discrete) numerical
applications, we use $\delta = 0$ in the measure~$\mu_\delta$ the trace
estimator vectors~$z_i$ are sampled from (see section~\ref{sec:form1}).

To enforce the constraint $\ww \in \SS$, i.e., $\sqrt{w_1^2 + w_2^2}=1$ in this
case, we
parameterize the weight vector as $(w_1, \pm \sqrt{1-w_1^2})$. The
parameter~$w_1$, alone, controls the combination of both sources. Moreover, the
weight vectors $(w_1, -\sqrt{1-w_1^2})$ and $(-w_1, \sqrt{1-w_1^2})$ lead to the
same reconstruction, such that it suffices to consider the encoding
weights~$(w_1, \sqrt{1-w_1^2})$ for $w_1 \in [-1,1]$.  

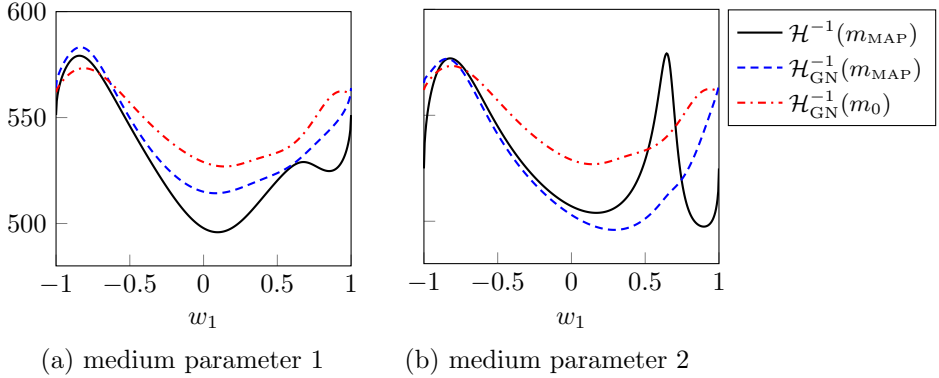
\begin{figure}
\centering
\begin{tabular}{cc}
\begin{tikzpicture}
\pgfplotstableread{a220--1123-10-10-11222-11101101-7-401.out}\mytable
\begin{axis}[width=\tikzw\textwidth, scale only axis,
	compat=newest,
	xlabel style={align=center},
	xlabel = $w_1$,
	xmin=-1.0, xmax=1.0, xtick pos=left,
	ylabel style={align=center},
	ylabel shift=-.1in,
	ytick pos=bottom,
	ymin=480, ymax=600, ylabel shift=-.05in,
    legend style={font=\small,nodes=right},
    legend pos = north west]
\addplot [color=black, thick]
table[x=w1, y=full] {\mytable};
\addplot [color=blue, thick, densely dashed]
table[x=w1, y=GN] {\mytable};
\addplot [color=red, thick, dashdotted]
table[x=w1, y=GNlin] {\mytable};
\end{axis}
\end{tikzpicture}
&
\begin{tikzpicture}
\pgfplotstableread{a220--2123-10-10-11222-11101101-7-401.out}\mytable
\begin{axis}[width=\tikzw\textwidth, scale only axis,
	xlabel style={align=center},
	xlabel = $w_1$,
	xmin=-1.0, xmax=1.0, xtick pos=left,
	ytick pos=bottom, ymin=480, ymax=600,
	yticklabel=\empty,
    legend style={font=\small,nodes=right},
    legend pos = outer north east]
\addplot [color=black, thick]
table[x=w1, y=full] {\mytable};
\addlegendentry{$\Hess^{-1}(m_\map)$}
\addplot [color=blue, thick, densely dashed]
table[x=w1, y=GN] {\mytable};
\addlegendentry{$\Hess^{-1}_{\text{GN}}(m_\map)$}
\addplot [color=red, thick, dashdotted]
table[x=w1, y=GNlin] {\mytable};
\addlegendentry{$\Hess^{-1}_{\text{GN}}(m_0)$}
\end{axis}
\end{tikzpicture}
\\
(a) medium parameter 1 & \hspace{-1.3in} (b) medium parameter 2
\end{tabular}
\caption{Plots of tr$(\Hess^{-1})$ with $\Hess^{-1}(m_\map(w_1))$,
$\Hess^{-1}_\text{GN}(m_\map(w_1))$ and $\Hess^{-1}_\text{GN}(m_0)$ for both
target media. $m_0 \equiv 1$, same as the background value for the medium
parameter.  }
\label{fig:1Dplots}
\end{figure}

In figure~\ref{fig:1Dplots}, we plot the three objective
functions~\eqref{eq:Wselgninf}--\eqref{eq:Wselgnlininf} from
section~\ref{sec:form1}. For each $w_1 \in [-1,1]$, the
Gauss--Newton~\eqref{eq:Wselgninf} and Laplace~\eqref{eq:Wsellaplaceinf}
formulations are evaluated at the MAP~point, $m_\map(w_1)$, corresponding to the
encoding weight~$(w_1, \sqrt{1-w_1^2})$; in other words, the Hessian for these
two criteria is evaluated at a medium parameter~$m_\map(w_1)$ that varies with
the weight~$w_1$.  For formulation~\eqref{eq:Wselgnlininf}, we choose $m_0$ to
be a constant value equal to the background medium, i.e., $m_0\equiv 1$.  
We observe that the result for the Gauss--Newton formulation~\eqref{eq:Wselgninf}
differs from the Laplace approximation~\eqref{eq:Wsellaplaceinf}.  In addition,
we clearly observe that each formulation contains local minima. 

\paragraph{Robustness of the Gauss--Newton formulation~\eqref{eq:Wselgninf}}

Since the computation of the MAP~point $m_\map(w_1)$ is a computationally
intensive task for large-scale problems, it might be useful to  
solve the optimization~\eqref{eq:Wselinf} without having
to recompute the exact MAP~point for each iterate of the weights. The Laplace
formulation~\eqref{eq:Wsellaplaceinf} is based on the full Hessian which is
guaranteed to be positive definite only in a neighbourhood of the MAP~point.
The Gauss--Newton approximation, however, is always positive definite and we
observe numerically that it preserves relevant information about the objective
function, even far away from the MAP~point.  In figure~\ref{fig:GN}, we plot the
objective function~\eqref{eq:Wselgninf}, for all values of $w_1 \in [-1,1]$, 
for different (fixed) medium parameters~$\bar{m}_s$ ranging from the background
medium,~$m_0 \equiv 1$, to the MAP~point~$m^\sharp$ computed using both sources
independently (for medium parameter~2).  The sources are located at
the points $(0,0.1)$ and~$(0,1.1)$.  That is, we define 
\[   \bar{m}_s = (1-s) m_0 + s\, m^\sharp. \] 

\newcommand{\mywsq}{0.15}
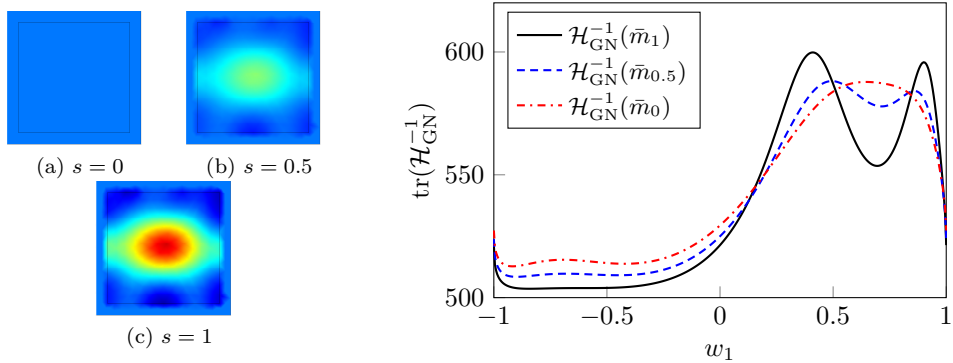
\begin{figure}
\centering
\begin{tabular}[c]{cl}
\begin{tabular}[c]{cc}
\multicolumn{2}{c}{\vspace{-.4in}} \\
\subfloat[$s=0$]{\includegraphics[width=\mywsq\textwidth, trim=12 10 50 10,
clip=true]{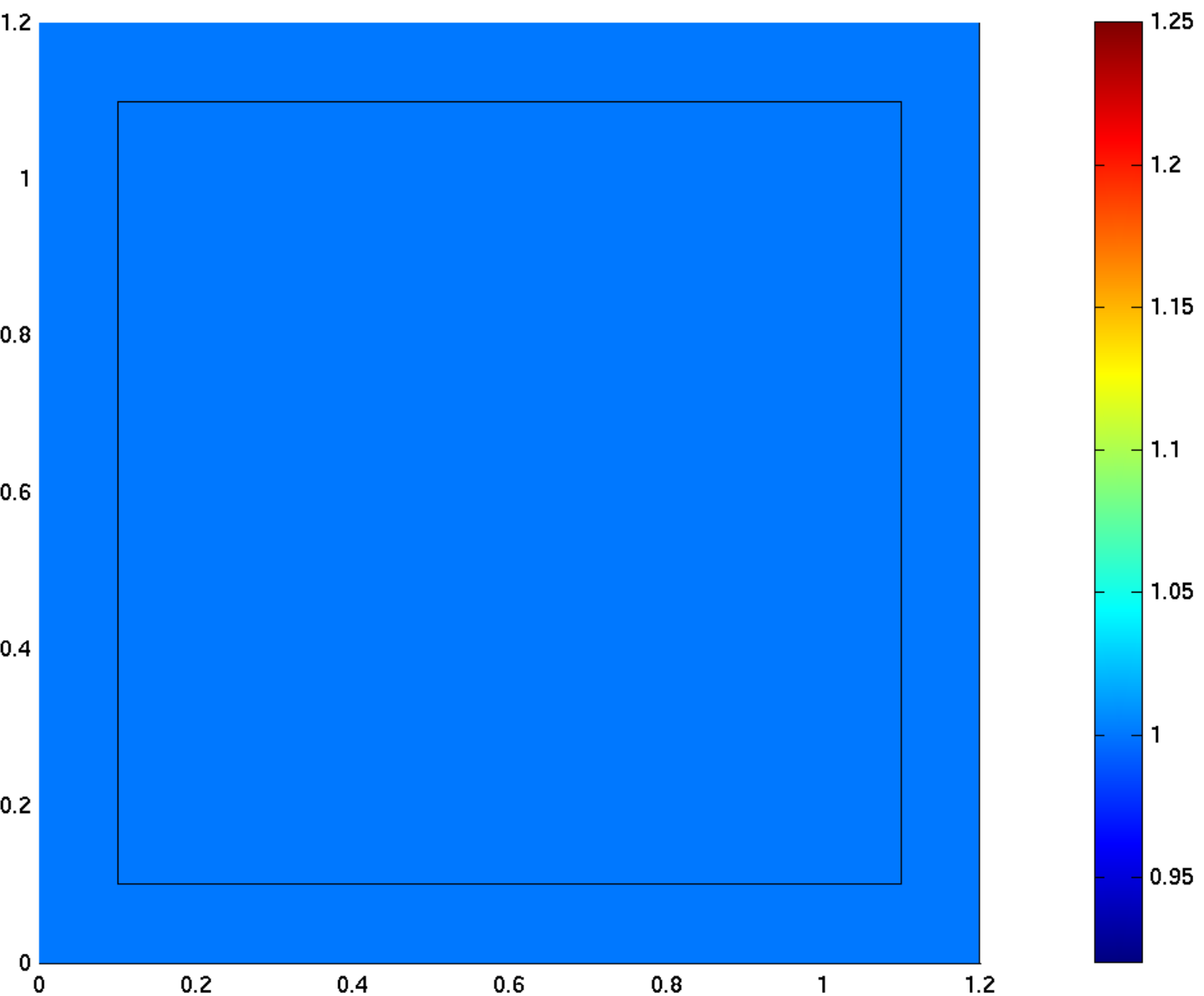}} &
\subfloat[$s=0.5$]{\includegraphics[width=\mywsq\textwidth, trim=12 10 50 10,
clip=true]{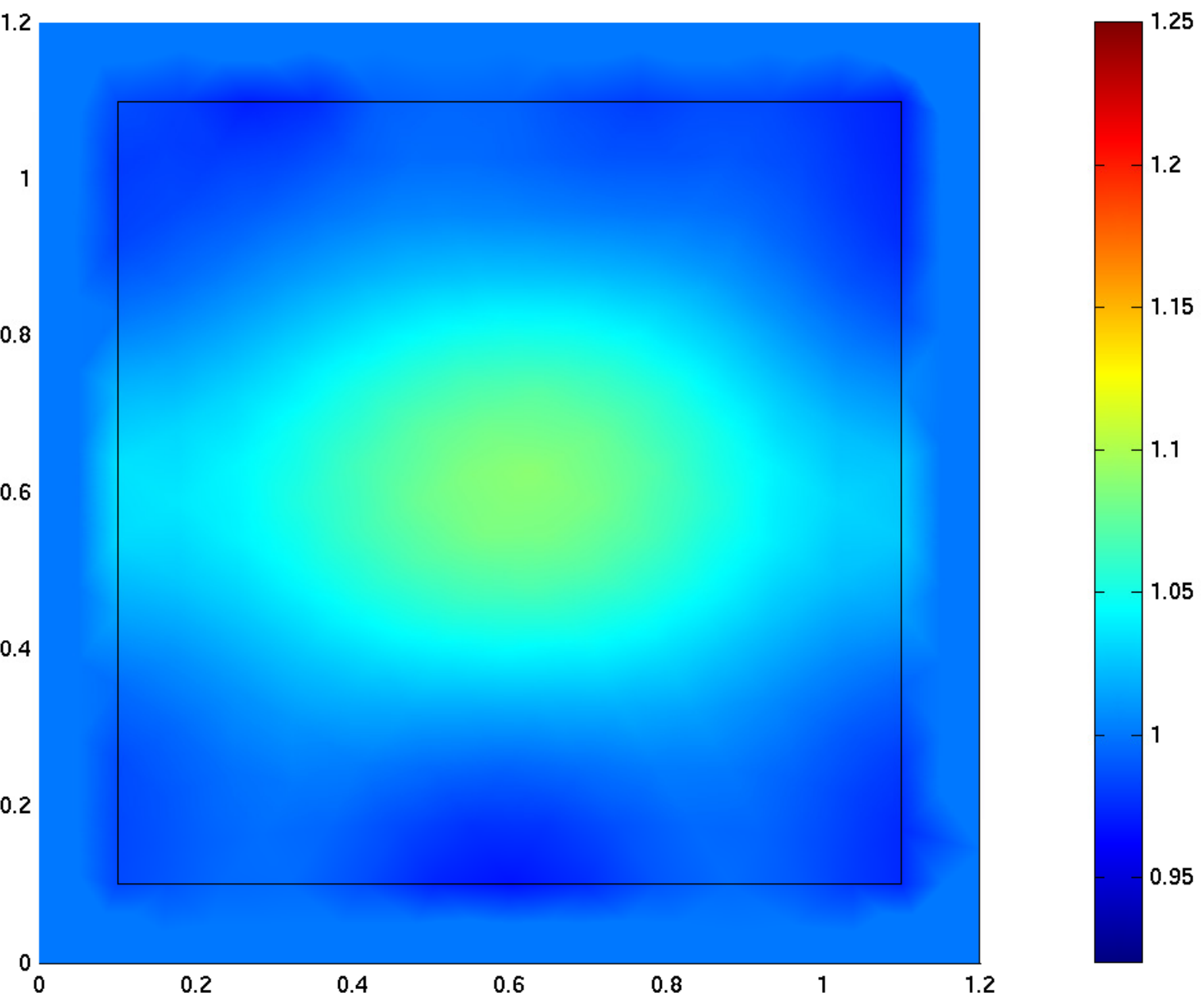}} \\
\multicolumn{2}{c}{\vspace{-.3in}} \\
\multicolumn{2}{c}{\subfloat[$s=1$]{\includegraphics[width=\mywsq\textwidth, trim=12 10 50 10,
clip=true]{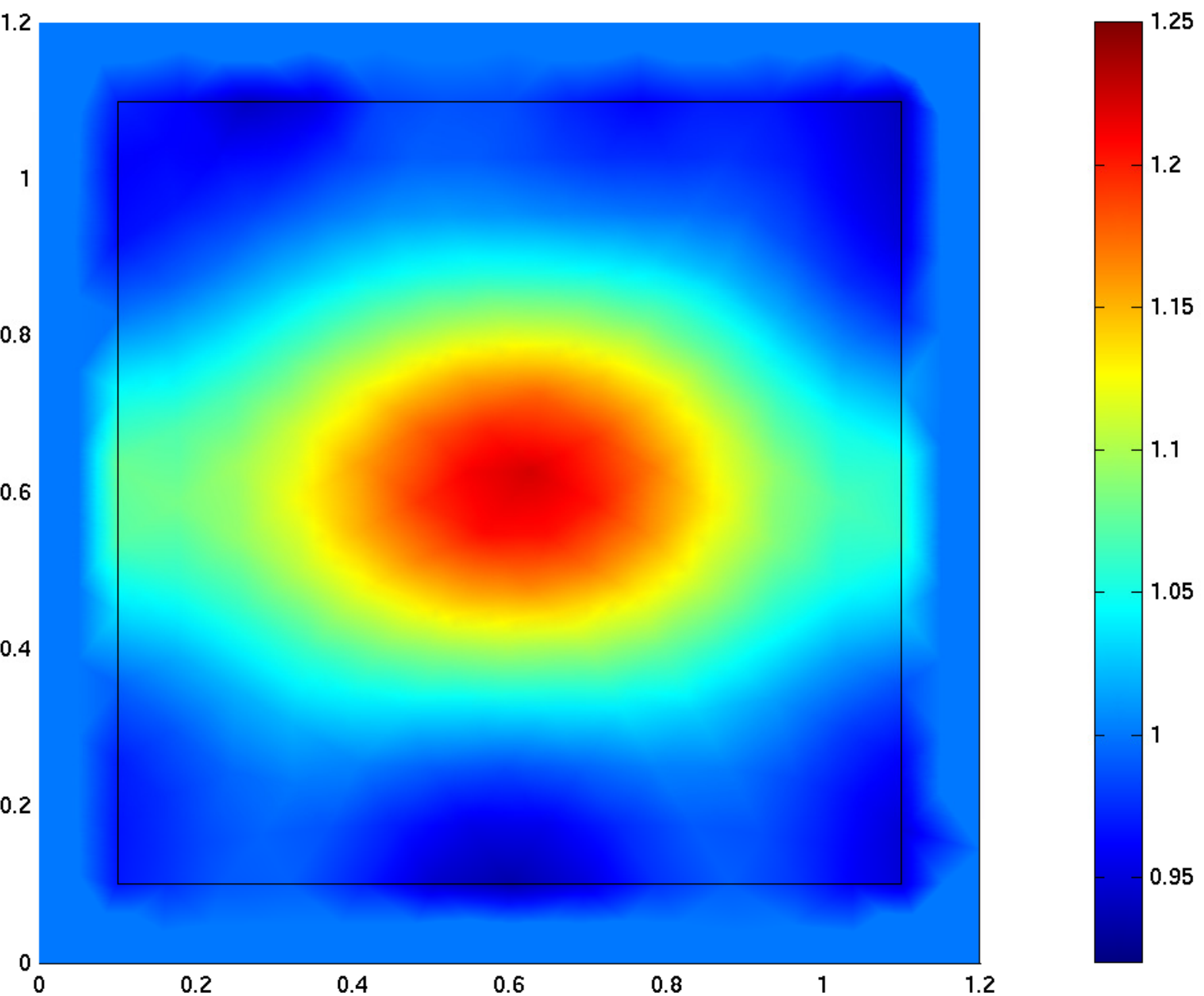}}}
\end{tabular}
&
\begin{tabular}[c]{c}
\begin{tikzpicture}
\begin{axis}[width=.46\textwidth, height=0.3\textwidth,scale only axis,
	compat=newest,
	xlabel = $w_1$, xmin=-1.0, xmax=1.0, xtick pos=left,
	ylabel = $\text{tr}(\Hess_\text{GN}^{-1})$, 
	ylabel shift=-.1in,
	ytick pos=bottom,
	ymin=500, ymax=620, ylabel shift=-.05in,
    legend style={font=\small,nodes=right},
    legend pos = north west]
\addplot [color=black, thick]
table[x=w1, y=GN]
{220--2123-10-10-11222-11101101-7-401-iter20.out};
\addlegendentry{$\Hess^{-1}_\text{GN}( \bar{m}_1)$}
\addplot [color=blue, thick, densely dashed]
table[x=w1, y=GN]
{220--2123-10-10-11222-11101101-7-401-iter10.out};
\addlegendentry{$\Hess^{-1}_\text{GN}( \bar{m}_{0.5})$}
\addplot [color=red, thick, dashdotted]
table[x=w1, y=GN]
{220--2123-10-10-11222-11101101-7-401-iter1.out};
\addlegendentry{$\Hess^{-1}_\text{GN}( \bar{m}_0)$}
\end{axis}
\end{tikzpicture} 
\end{tabular}
\end{tabular}
\caption{Plots of objective function $\obj_0$~\eqref{eq:Wselgnlininf} for
weights~$w_1 \in [-1,1]$ (right), at medium $\bar{m}_s$, with $s = 0, \, 0.5, \, 1$
(left).  Here $m_0 \equiv 1$ (the background medium).
}
\label{fig:GN}
\end{figure}

It appears that the medium parameter needs to include the main features of the
target medium sufficiently accurately~($s > 0.5$) to match the main features of
the exact trace of the posterior covariance; this can be seen from
the behavior of $\text{tr}\big(\Hess^{-1}_\text{GN}(w_1, \bar{m}_s)\big)$ in the
interval~$w_1 \in [0.2,1.0]$.

\paragraph{The effect of trace estimation}

\begin{figure}
\centering
\begin{tikzpicture}
\pgfplotstableread{220--2123-10-10-11222-11101101-7-401-trestim.out}\trestim
\begin{axis}[width=.35\textwidth, scale only axis,
	compat=newest,
	xmin=-1.0, xmax=1.0, xtick pos=left,
	ylabel style={align=center},
	ylabel = $\text{tr}(\Hess^{-1})$, 
	ylabel shift=-.1in,
	ytick pos=bottom,
	ymin=300, ymax=600, ylabel shift=-.05in,
    legend style={font=\small,nodes=right},
    legend pos = outer north east]
\addplot [color=gray, ultra thick]
table[x=w1, y=full] {\trestim};
\addlegendentry{exact}
\addplot [color=black, thick]
table[x=w1, y=Trest30] {\trestim};
\addlegendentry{$n_{tr}=30$}
\addplot [color=blue, thick, densely dashed]
table[x=w1, y=Trest10] {\trestim};
\addlegendentry{$n_{tr}=10$}
\addplot [color=red, thick, dashdotted]
table[x=w1, y=Trest1] {\trestim};
\addlegendentry{$n_{tr}=1$}
\end{axis}
\end{tikzpicture} 
\caption{Plots of the objective function in~\eqref{eq:Wclaplaceinf} when
the trace of the posterior covariance is computed exactly or with a
trace estimator ($n_{tr}=1, \, 10, \, 30$).
For each $n_{tr}$, we used a fixed realization of the trace estimator vectors.}
\label{fig:trestimfig}
\end{figure}
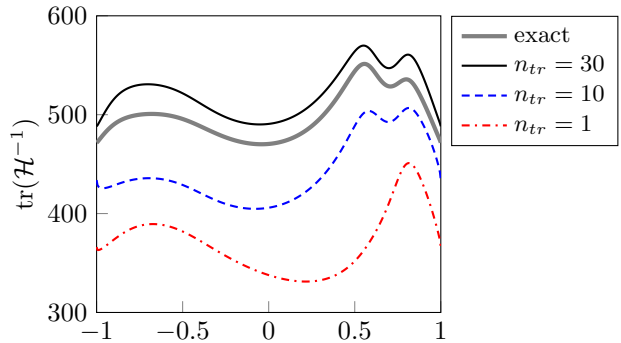

When computing A-optimal encoding weights, 
one only needs the local minima of the trace to be well characterized. 
We show in figure~\ref{fig:trestimfig} that trace estimation does indeed affect
the shape of the objective function in the formulations of the A-optimal
encoding weights~\eqref{eq:Wclaplaceinf}.  However, in our example, the
objective function using a trace estimation preserves the local minima  of the
objective function using an exact trace when a sufficient number of trace
estimator vectors are used.

\subsection{A-optimal encoding weights in higher dimensional weight spaces}

We now consider a problem with 10~sources (i.e., $N_s=10$).  
\begin{figure}[h]
\centering
\begin{tikzpicture}
\node[anchor=south west, inner sep=0] (image1) at (0,0)
{\includegraphics[height=0.25\textwidth, trim=50 20 90 10,
clip=true]{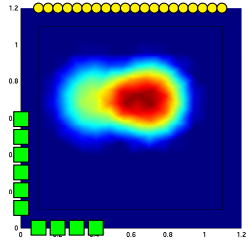}};
\begin{scope}[x={(image1.south east)},y={(image1.north west)}]
\foreach \x in {0,1,...,5} {
\filldraw[fill=green,draw=black] (0.05,0.47-.072*\x) rectangle (0.11,0.53-.072*\x);}
\foreach \x in {0,1,...,3} {
\filldraw[fill=green,draw=black] (0.12+0.077*\x,0.03) rectangle (0.18+0.077*\x,0.09);}
\foreach \x in {0,1,...,19} {
\filldraw[fill=yellow,draw=black] (0.15+0.03895*\x,.95) circle (0.02);}
\end{scope} 
\end{tikzpicture}
\caption{Target medium parameter and locations of the 10~sources (green
squares), and receivers (yellow circles).}
\label{fig:dom10src} \end{figure}
Here, we focus on qualitative properties of the A-optimal source encoding
weights by performing statistical tests, in which we study how successful
A-optimal encoding weights are in reducing posterior variance and relative
medium misfit compared to encoding weights sampled from the uniform spherical
distribution.  We also compared with random weights sampled, then re-scaled,
from the Rademacher distribution (see section~\ref{sec:medparam}). Since the
results we obtained were not statistically different from the results presented
in this section using random weights sampled from the uniform spherical
distribution, we decided to omit these results.  Throughout this section, the
relative medium misfit is taken to be the relative $L^2$-error between the
reconstruction of interest and the reconstruction obtained using all 10~sources
independently.  The penalty parameter was empirically selected to
be~$\lambda=10^3$.

\begin{figure}[htb]
\centering
\includegraphics[trim=100 410 110 90, clip=true]{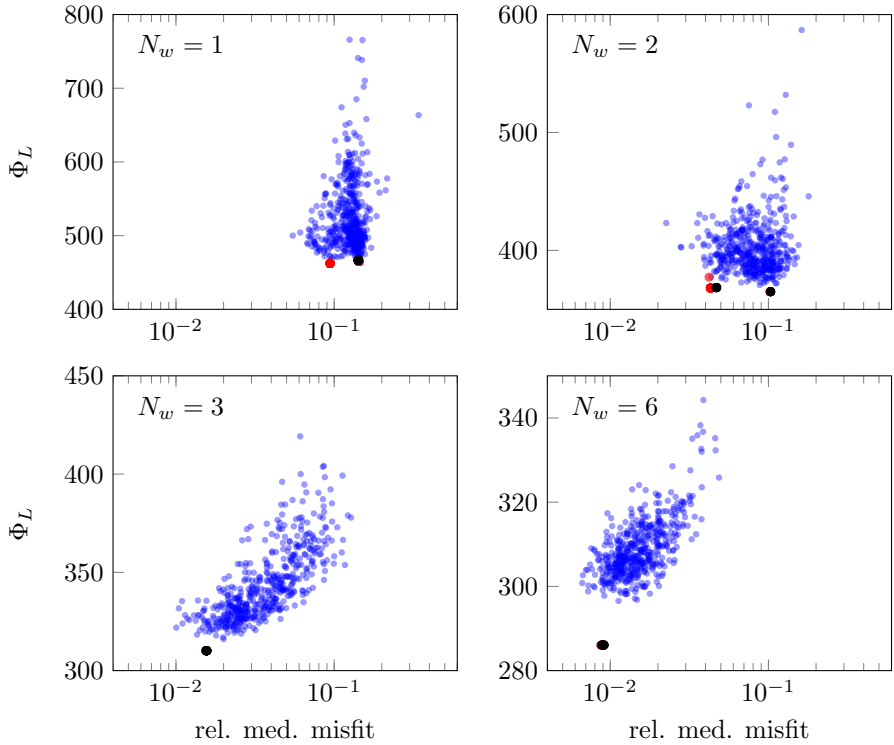}
\caption{Plot of~$\obj_L(\ww)$~\eqref{eq:Wsellaplaceinf} against relative medium
misfit ($N_s=10$ and $N_w=1, \, 2, \, 3, \, 6$) for reconstructions using random encoding
sources sampled from the uniform spherical distribution (blue) or A-optimal
encoding weights computed with formulation~\eqref{eq:Wcgninf}~(black)
and~\eqref{eq:Wclaplaceinf}~(red).  Target model~2 with source configuration as
shown in figure~\ref{fig:dom10src}.  Sample size = 500, $n_{tr}=30$.}
\label{fig:clouda10Ns} 
\end{figure}

We show the results in figure~\ref{fig:clouda10Ns}.  Each plot shows, for
different number of encoded sources~($N_w=1$, $2$, $3$ and $6$), the objective
function~$\obj_L(\ww)$ defined in~\eqref{eq:Wsellaplaceinf} against the relative
medium misfit of the reconstruction, which is an indication for the quality of
the reconstruction.  Each reconstruction is indicated by a translucent dot; a
darker shade indicates a higher concentration of reconstructions in that part of
the plot.  This shows the variation in the quality of the reconstruction.  The
blue dots correspond to reconstructions that use random encoding weights sampled
from the uniform spherical distribution.  The red dots indicate A-optimal
encoding weights based on the Laplace formulation~\eqref{eq:Wclaplaceinf}.  The
reconstructions marked with black dots use A-optimal encoding weights based on
the Gauss--Newton formulation~\eqref{eq:Wcgninf}.  In order to detect potential
local minima, the A-optimal encoding weights are re-computed several times,
starting from different initial conditions.

\begin{figure}[!h]\centering
\begin{tabular}{c @{\hspace{-.01in}} c @{\hspace{-.01in}} c @{\hspace{-.01in}} c}
\includegraphics[width=0.24\textwidth,
trim=5.4cm 9.25cm 6.35cm 8.9cm,
clip=true]{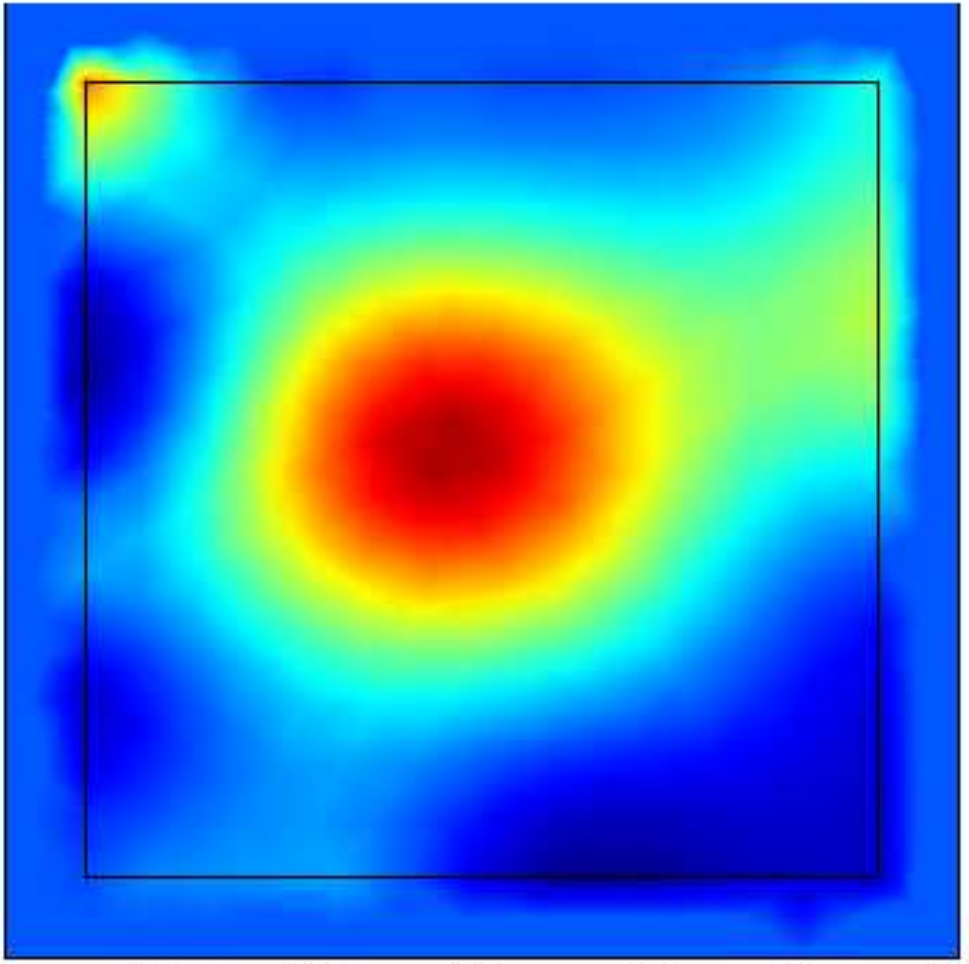} \label{fig:reconstruct10src}
& 
\includegraphics[width=0.24\textwidth,
trim=5.4cm 9.25cm 6.35cm 8.9cm,
clip=true]{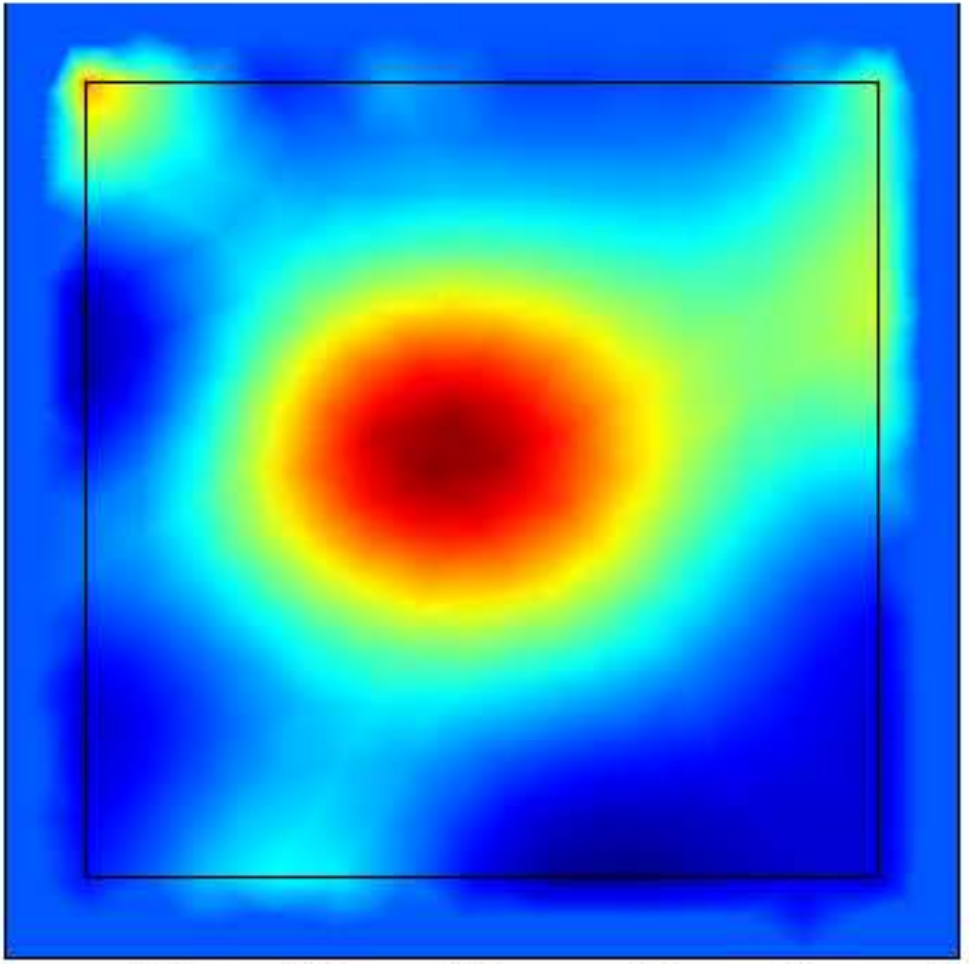} \label{fig:reconstructAopt}
&
\includegraphics[width=0.24\textwidth,
trim=5.4cm 9.25cm 6.35cm 8.9cm,
clip=true]{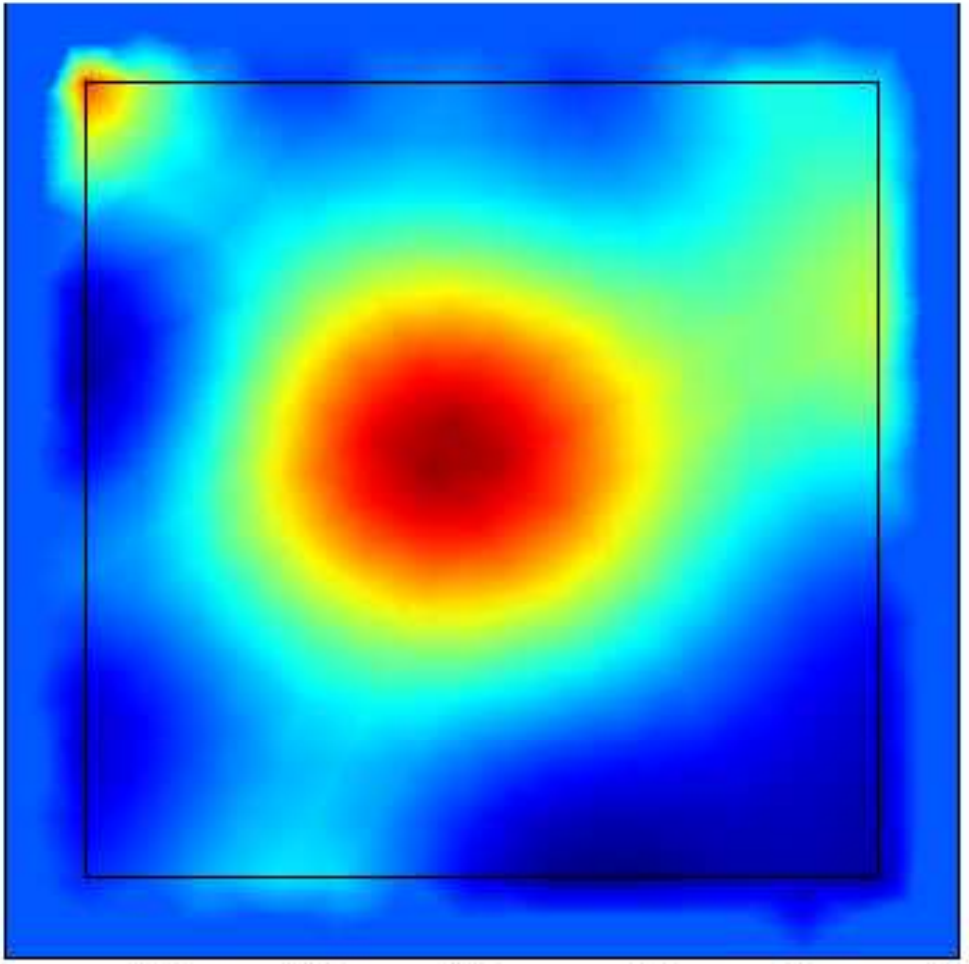} \label{fig:reconstructrndmin}
&
\includegraphics[width=0.24\textwidth,
trim=5.4cm 9.25cm 6.35cm 8.9cm,
clip=true]{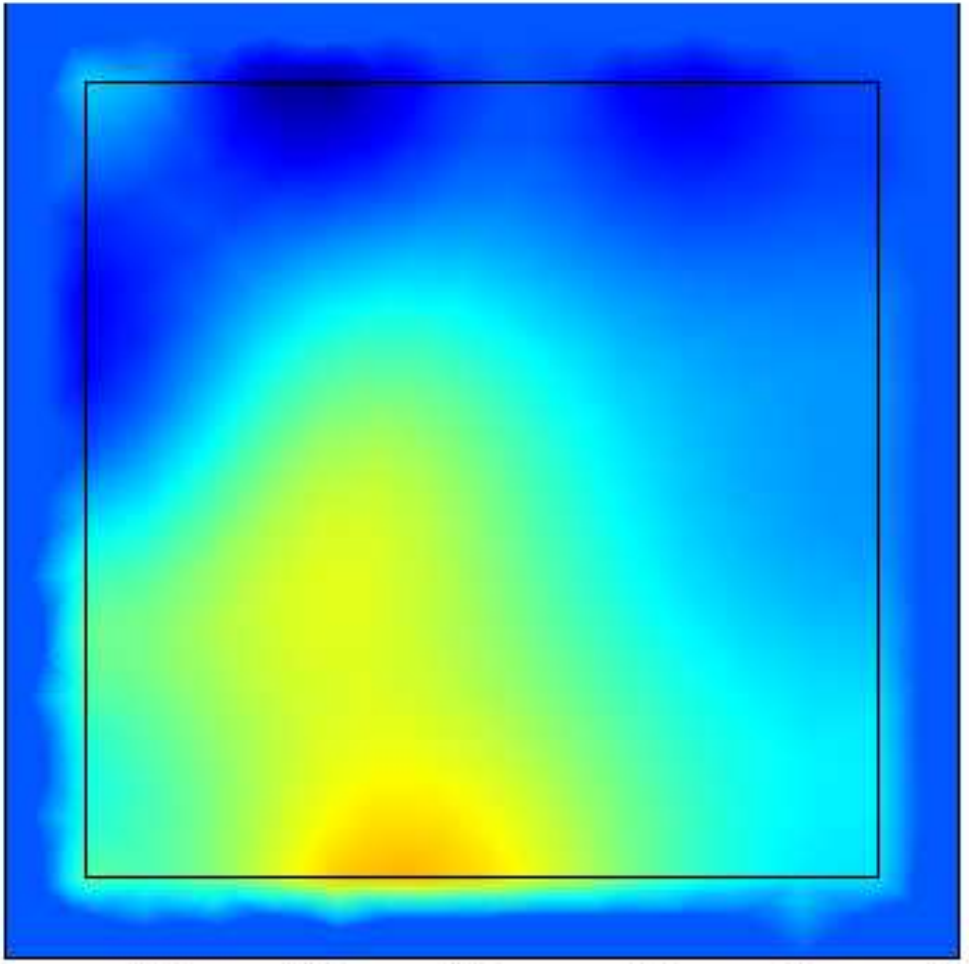} \label{fig:reconstructrndmax}\\
(a) & (b) & (c) & (d) 
\end{tabular}
\caption{Four examples of reconstructions using different number of sources,
with target parameter~2:
(a) 10 independent sources; 
(b) 3 A-optimally encoded sources;
(c) 3 randomly encoded sources;
(d) 3 other randomly encoded sources.
\label{fig:reconstruct}
}
\end{figure}

Notice that with one encoded source, A-optimal encoding weights do not provide a
clear advantage over random weights. The overall distribution of random weights
does not indicate a strong connection between the trace of the posterior
covariance~\eqref{eq:Wsellaplaceinf} and the relative medium misfit.  On the
other hand, the A-optimal encoding weights outperform the random weights (on
average), when sufficiently many encoding weights are used (see in particular
$N_w=2$ and~$3$ in figure~\ref{fig:clouda10Ns}).  In that case, the random
weights appear to indicate a linear correlation between our objective function
and the relative medium misfit, which translates into the best reconstruction
being also the one with smallest trace of the posterior covariance.  Overall,
these results suggest the existence of a threshold, in the number of encoding
sources, above which optimal weights provide improvement in both variance and
medium misfit over random encoding weights.  
Moreover, based on these results, there does not appear to be a clear advantage
in using the Laplace approximation~\eqref{eq:Wclaplaceinf} over the Gauss--Newton
approximation~\eqref{eq:Wcgninf}, provided sufficiently many encoded sources are
used. In the last row of figure~\ref{fig:clouda10Ns}, optimal weights computed
with both formulations provide similar results, although the actual values of
the weights do not necessarily agree.

In addition, we provide a comparison of the reconstructions computed using all
sources independently (figure~\ref{fig:reconstruct}a), using three~A-optimally
encoded sources (figure~\ref{fig:reconstruct}b), and two examples of
reconstructions computed using three randomly encoded sources: one resulting in
a good reconstruction (figure~\ref{fig:reconstruct}c), and one resulting in  a
poor reconstruction (figure~\ref{fig:reconstruct}d).  There is virtually no
difference between the reconstructions computed using all 10~sources and using
three A-optimally encoded sources.  On the other hand, using random encoding
weights drawn from the same distribution may lead to good or poor
reconstructions, as is shown in figures~\ref{fig:reconstruct}c, d. This is
consistent with the results in figure~\ref{fig:clouda10Ns} 
(bottom left), where the blue dots show large variations in terms of relative
medium misfit.

\paragraph{Variability of the A-optimal encoding weights}

The A-optimal encoding weight formulation introduced in
section~\ref{sec:formAweights} relies on a fixed realization of the trace
estimator vectors.  Note that the A-optimal encoding weights are solutions to a highly
nonlinear optimization problem that in general exhibits local minima.  However,
we show numerically that, provided sufficiently many encoding weights are chosen
and a large enough number of trace estimator vectors are used, the computation of the
A-optimal encoding weights is stable with respect to trace estimation.  
\begin{figure}[h]
\centering
\includegraphics[trim=110 510 110 90, clip=true]{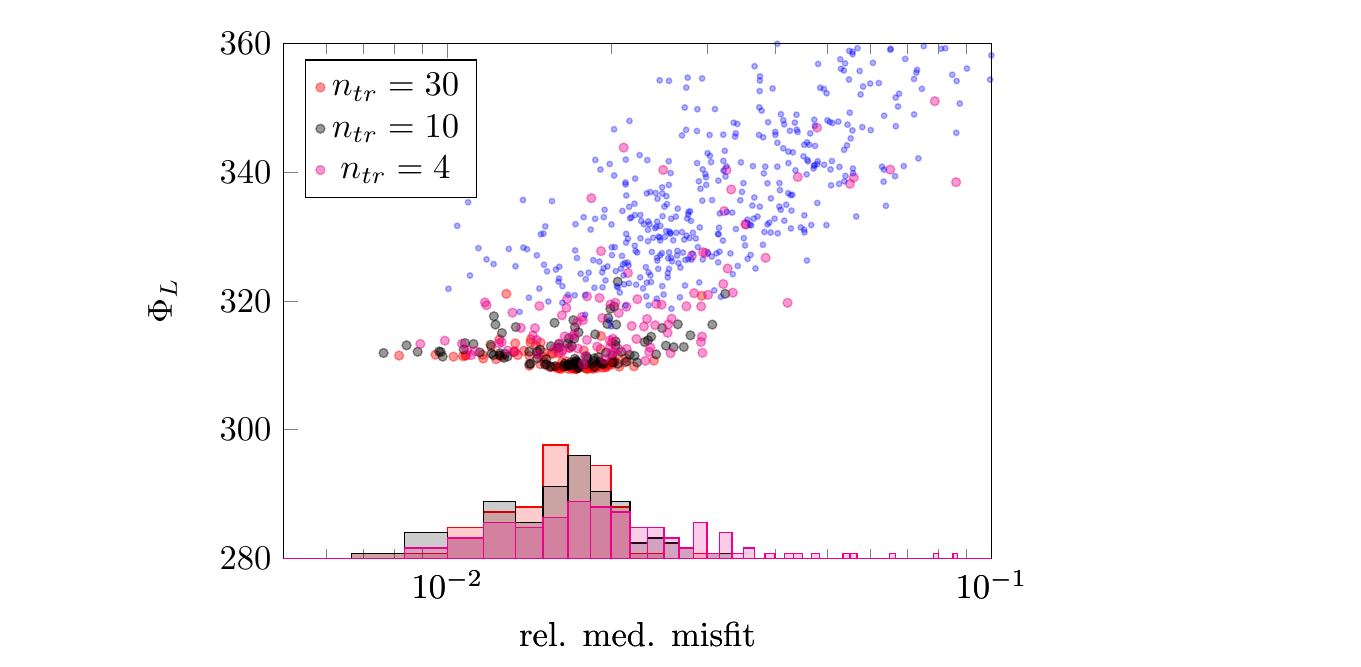}
\caption{
Variability of the A-optimal weights for different numbers of trace
estimator vectors, $n_{tr} = 30$ (red), $10$~(black) and $4$~(magenta).
A-optimal encoding weights are computed with formulation~\eqref{eq:Wclaplaceinf}
($N_s=10$ and $N_w=3$), using different realizations of the trace estimator
vectors and
different initial guess of the weights. Sample size = 100.}
\label{fig:trace} \end{figure}
In figure~\ref{fig:trace}, we show 100~results obtained with Laplace A-optimal
encoding weights~\eqref{eq:Wclaplaceinf}, in the case of 3~encoded sources, with
different numbers of trace estimator vectors ($n_{tr} = 4, \, 10, \, 30$).  Each computation
uses different realizations of the trace estimator vectors, and different initial guess
of the weights.

We observe that with $n_{tr}=10$ and~$30$ the computations of the A-optimal
encoding weights provide similar results. On the other hand, the use of
$4$~trace estimator vectors leads to a much wider range in the quality of the
results, both in terms of relative medium misfit and trace of the posterior
covariance.

\subsection{Remarks on the Gauss--Newton formulation}

Here, we discuss the justification for and advantages of using the Gauss--Newton
formulation for finding A-optimal encoding weights.  In many important
situations, the Gauss--Newton formulation appears accurate enough to compute the
A-optimal encoding weights.  The Gauss--Newton approximation to the Hessian is
most accurate when the data misfit residual is small at the solution of the
inverse problem. This is the case, for instance, when the noise level in the
observations is low.  In our numerical experiments we observed that, provided
sufficiently many encoded sources are used, the Gauss--Newton formulation
represents a sufficiently accurate approximation to the Laplace formulation for
the purpose of computing A-optimal encoding weights.  

The Gauss--Newton formulation holds strong promises to reduce the computational
cost of the A-optimal encoding weights.  The data-misfit part of the
Gauss--Newton Hessian is guaranteed to be positive semi-definite at any
evaluation point, and hence the Gauss--Newton Hessian is positive definite. This
allows two main improvements to the computations of the A-optimal weights.
First, and as detailed in section~\ref{sec:computational_cost}, one can
incorporate a low-rank approximation of the Gauss--Newton Hessian to reduce the
computational cost.  The magnitude of that reduction is problem-dependent, but
will be most noticeable when large numbers of trace estimator vectors are required.

Another advantage of the positive definiteness of the Gauss--Newton Hessian is
that the objective function~\eqref{eq:Wselgninf} of the Gauss--Newton formulation
does not have to be evaluated in a small neighbourhood of the MAP~point for the
objective function to make sense.  This could allow one, for instance, to solve
the MAP~point inexactly when the A-optimal objective function is far from its
minimum, which would reduce the overall computational cost.  In
section~\ref{sec:1dspace}, we studied how the objective function varies with the
evaluation point~$\bar{m}_s$ (figure~\ref{fig:GN}), and observed that the
objective function tends to maintain similar local minima away from the MAP
point. 

Finally, we want to point out that in certain situations, the full Hessian may
not be available, may be too complicated to derive, or too expensive to compute,
rendering the Laplace formulation inadequate.  This can be the case for inverse
problems with highly nonlinear forward problems. 

\section{Conclusion}
\label{sec:conc}

We have developed a method for the computation of A-optimal encoding weights
aiming at large-scale non-linear inverse problems.  As we show numerically,
reconstructions obtained using A-optimal encoding weights not only minimize the
average of the posterior variance, but consistently outperform random encoding
weights in terms of the quality of the reconstructions.
While in this work, we relied on quasi-Newton methods for solving the
optimization problem for A-optimal encoding weights, we will explore the
derivation and implementation of a Newton solver for this optimization problem
in future work.  We point out that, thanks to the optimize-then-discretize
approach we adopted, the derivation of the analytical expression for the action
of the Hessian in a direction is possible with little more effort than what was
required to get the gradient.

We introduced two formulations for the computation of the A-optimal encoding
weights, namely the Gauss--Newton formulation~\eqref{eq:Wcgninf} and the Laplace
formulation~\eqref{eq:Wclaplaceinf}. Although the Gauss--Newton formulation
represents an approximation to the Laplace formulation, it holds several
advantageous features from computational point of view.

We note that computing A-optimal encoding weights can entail a significant
computational effort.  However, the method can be attractive for real-time
monitoring applications where one needs to solve an inverse problem repeatedly
over time.  In this case, one first computes the A-optimal encoding weights
offline, and then can use those weights to solve the inverse problem repeatedly
at a fraction of the original cost.  An example for such an application is the
monitoring of an oil reservoir, where seismic or electro-magnetic inverse
problems are solved repeatedly to characterize the evolution of the reservoir
properties over time.

\section*{References}
\bibliographystyle{unsrt}

\appendix
\section{Gradient of the optimization formulation~\eqref{eq:wkform}}
\label{sec:outer}

We detail the derivation of the gradient of the Laplace formulation of the
A-optimal weights in the case of the Helmholtz inverse problem, as defined
in~\eqref{eq:wkform}.  In that formulation, we enforce the PDE constraints
weakly using Lagrange multipliers.  Therefore, we need to introduce adjoint
variables that are indicated with a star superscript, e.g., $m^*$ is the adjoint
variable for~$m$.  Following the formal Lagrangian approach~\cite{Troltzsch10},
we define the Lagrangian~$\Lagr$,
\begin{align} 
& \hspace{.2in} \Lagr (\ww,m,\{u_i\},\{p_i\},\{v_{i,k}\},\{q_{i,k}\},\{y_{k}\},
m^*,\{u_i^*\},\{p_i^*\},\{v_{i,k}^*\},\{q_{i,k}^*\},\{y_{k}^*\}) = \nonumber \\
 & \frac1{n_{tr}} \sum_{k=1}^{n_{tr}} \langle y_k, z_k \rangle + \nonumber \\
& \frac1{n_{tr}N_w} \sum_{k=1}^{n_{tr}} \sum_{i=1}^{N_w} \left[ \langle \nabla
v_{i,k}, \nabla v_{i,k}^* \rangle - \kappa^2 \langle m v_{i,k}, v_{i,k}^*
\rangle - \kappa^2 \langle u_i y_k, v_{i,k}^* \rangle  \right] \nonumber \\
& + \frac1{n_{tr}N_w} \sum_{k=1}^{n_{tr}} \sum_{i=1}^{N_w} \bigg[ \langle \nabla
q_{i,k}^*, \nabla q_{i,k} \rangle - \kappa^2 \langle q_{i,k}^*, m q_{i,k}
\rangle - \kappa^2 \langle q_{i,k}^*, p_i y_k \rangle + \langle B q_{i,k}^*, B
v_{i,k} \rangle_{\GG^{-1}_\text{noise}} \bigg] \nonumber \\
& + \frac1{n_{tr}} \sum_{k=1}^{n_{tr}} \left[ \langle y_k, y_k^* \rangle_\E
 - \frac1{N_w} \sum_{i=1}^{N_w} \kappa^2 \left( \langle v_{i,k} p_i, y_k^*
   \rangle + \langle  u_i q_{i,k}, y_k^* \rangle \right) 
- \langle z_k, y_k^* \rangle \right] + \nonumber \displaybreak[0] \\
& \frac1{N_w} \sum_{i=1}^{N_w} \left[ \langle \nabla u_i , \nabla u_i^* \rangle 
- \kappa^2 \langle m u_i , u_i^* \rangle - \langle f(\ww^i), u_i^* \rangle
  \right] \nonumber \\
& + \frac1{N_w} \sum_{i=1}^{N_w} \left[ \langle \nabla p_i^*, \nabla p_i \rangle
- \kappa^2 \langle p_i^*, m p_i \rangle + \langle B p_i^*,  B u_i - \dd(\ww^i)
  \rangle_{\GG^{-1}_\text{noise}} \right] \nonumber \\
& + \langle m -m_0 , m^* \rangle_\E 
- \frac1{N_w} \sum_{i=1}^{N_w} \kappa^2 \langle u_i p_i, m^* \rangle .
\label{eq:LagrO}
\end{align} 
The gradient is then given by $\delta_\ww \Lagr = \left[ \delta_{\ww^1} \Lagr,
\delta_{\ww^2} \Lagr, \dots, \delta_{\ww^{N_w}} \Lagr \right]^T$, where for
any $i=1, \ldots, N_w$,
\begin{equation*}
\delta_{\ww^i} \Lagr =  - \frac1{N_w} 
\begin{bmatrix}
\langle f_1, u_i^* \rangle & + &
\langle B p_i^*, \dd_1 \rangle_{\GG^{-1}_\text{noise}} \\
\langle f_2, u_i^* \rangle & + &
\langle B p_i^*, \dd_2 \rangle_{\GG^{-1}_\text{noise}} \\
& \vdots & \\
\langle f_{N_s}, u_i^* \rangle & + &
\langle B p_i^*, \dd_{N_s} \rangle_{\GG^{-1}_\text{noise}} 
\end{bmatrix}. \end{equation*}

Before we specify the steps that lead to the evaluation of the
variables~$u_i^*$ and~$p_i^*$, we identify some important symmetries between
the state variables and their adjoints.  Indeed, for each $k=1, \ldots, n_{tr}$,
the variables~$(y_k \{v_{i,k}\}_i, \{q_{i,k}\}_i)$ solve a Hessian system similar
to~\eqref{eq:Hz}, and the corresponding adjoint variables~$(y_k^*
\{v^*_{i,k}\}_i, \{q^*_{i,k}\}_i)$ solve the system of
equations given (formally) by $\delta_{v_{ik}} \Lagr = \delta_{q_{ik}} \Lagr
= \delta_{y_k} \Lagr = 0$. While the former system of equations solve~$\Hess
y_k = z_k$, the latter solves~$\Hess y_k^* = - z_k$. This leads to the
symmetry relations 
\begin{equation} \label{eq:symm} y_k = - y_k^*, \ v_{ik} = -q_{ik}^*, \text{
and } q_{ik} = -v_{ik}^*, \end{equation}
for any $i=1, \ldots, N_w$ and $k=1, \ldots, n_{tr}$.

For any $i=1, \ldots, N_w$, the variable~$u_i^*$ (resp.~$p_i^*$) solves the
equation~$\delta_{u_i} \Lagr = 0$ (resp.~$\delta_{p_i} \Lagr = 0$).
That is, for any $\tilde{u} \in H^1(\dom)$, $u_i^*$ solves
\begin{multline*}
\langle \nabla u_i^*, \nabla \tilde{u}  \rangle 
- \kappa^2 \langle m u_i^*, \tilde{u} \rangle \\
- \kappa^2 \langle p_i m^*, \tilde{u} \rangle
+ \langle B p_i^*, B \tilde{u} \rangle_{\GG^{-1}_\text{noise}}
- \kappa^2 \frac1{n_{tr}} \sum_{k=1}^{n_{tr}} \left[ \langle y_k v_{i,k}^*, \tilde{u}
\rangle + \langle q_{i,k} y_k^*, \tilde{u} \rangle
\right] = 0.
\end{multline*} 
On the other hand, for any $\tilde{p} \in H^1(\dom)$, $p_i^*$ solves
\[ \langle \nabla p_i^*, \nabla \tilde{p} \rangle
- \kappa^2 \langle p_i^*, m \tilde{p} \rangle 
- \kappa^2 \langle u_i m^*, \tilde{p} \rangle 
- \kappa^2 \frac1{n_{tr}} \sum_{k=1}^{n_{tr}} \left[
\langle q_{i,k}^*y_k , \tilde{p} \rangle + \langle v_{i,k} y_k^*, \tilde{p}
\rangle \right] = 0. \]
Using~\eqref{eq:symm}, this reduces, for any $i=1, \ldots, N_w$,
to the system of equations
\begin{equation}  \label{eq:uipistar} \begin{aligned}
\langle \nabla u_i^*, \nabla \tilde{u}  \rangle 
- \kappa^2 \langle m u_i^*, \tilde{u} \rangle 
- \kappa^2 \langle p_i m^*, \tilde{u} \rangle
+ \langle B p_i^*, B \tilde{u} \rangle_{\GG^{-1}_\text{noise}}
+ \frac2{n_{tr}} \sum_{k=1}^{n_{tr}} \kappa^2 \langle y_k q_{i,k}, \tilde{u}
\rangle & = 0, \\
\langle \nabla p_i^*, \nabla \tilde{p} \rangle
- \kappa^2 \langle m p_i^*, \tilde{p} \rangle 
- \kappa^2 \langle u_i m^*, \tilde{p} \rangle 
+ \frac2{n_{tr}} \sum_{k=1}^{n_{tr}} \kappa^2 \langle v_{i,k} y_k , \tilde{p}
\rangle & = 0 .
\end{aligned} \end{equation}
Therefore, computation of the $u_i^*$'s and $p_i^*$'s requires knowledge of the
quantities $\{ u_i\}$, $\{p_i\}$, $m^*$, $\{v_{i,k}\}$, $\{q_{i,k}\}$ and $\{y_k\}$.
Variables~$\{ u_i\}$, $\{p_i\}$, $\{v_{i,k}\}$, $\{q_{i,k}\}$, and $\{y_k\}$ are
all evaluated during the computation of the objective functional~$1/n_{tr}
\sum_{k=1}^{n_{tr}} \langle y_k, z_k \rangle$, such that the only remaining
unknown quantity is~$m^*$. That variable is solution to the equation $\delta_m
\Lagr = 0$, that is, for any $\tilde{m} \in \E$, $m^*$ solves
\begin{multline*} 
\frac1{n_{tr}N_w} \sum_{k=1}^{n_{tr}} \sum_{i=1}^{N_w} \left[
- \kappa^2 \langle \tilde{m} v_{i,k}, v_{i,k}^* \rangle 
- \kappa^2 \langle q_{i,k}^*, \tilde{m} q_{i,k} \rangle \right] \\
+ \frac1{N_w} \sum_{i=1}^{N_w} \left[
- \kappa^2 \langle \tilde{m} u_i , u_i^* \rangle 
- \kappa^2 \langle p_i^*, \tilde{m} p_i \rangle \right]
+ \langle \tilde{m} , m^* \rangle_\E = 0.
\end{multline*}
Using~\eqref{eq:symm}, we simplify this equation to obtain
\[ \frac2{n_{tr}N_w} \sum_{k=1}^{n_{tr}} \sum_{i=1}^{N_w} 
\kappa^2 \langle v_{i,k} q_{i,k}, \tilde{m} \rangle 
- \frac1{N_w} \sum_{i=1}^{N_w} \kappa^2 \left[
 \langle u_i u_i^*, \tilde{m} \rangle 
+ \langle p_i^*p_i,\tilde{m} \rangle \right]
+ \langle m^*, \tilde{m} \rangle_\E = 0. \]
This equation can be grouped with the system of equations~\eqref{eq:uipistar} to
obtain the larger system
\begin{align*}
\langle \nabla p_i^*, \nabla \tilde{p} \rangle &
- \kappa^2 \langle m p_i^*, \tilde{p} \rangle 
- \kappa^2 \langle u_i m^*, \tilde{p} \rangle 
 = - \frac2{n_{tr}} \sum_{k=1}^{n_{tr}} \kappa^2 \langle v_{i,k} y_k , \tilde{p}
\rangle \\
\langle \nabla u_i^*, \nabla \tilde{u}  \rangle &
- \kappa^2 \langle m u_i^*, \tilde{u} \rangle 
- \kappa^2 \langle p_i m^*, \tilde{u} \rangle 
+ \langle  B p_i^*, B \tilde{u} \rangle_{\GG^{-1}_\text{noise}}
 = - \frac2{n_{tr}} \sum_{k=1}^{n_{tr}} \kappa^2 \langle y_k q_{i,k}, \tilde{u}
\rangle \\
\langle m^*, \tilde{m} \rangle_\E  &
- \frac1{N_w} \sum_{i=1}^{N_w} \kappa^2 \left[
 \langle u_i u_i^*, \tilde{m} \rangle 
+ \langle p_i^*p_i,\tilde{m} \rangle \right]  =
- \frac2{n_{tr}N_w} \sum_{k=1}^{n_{tr}} \sum_{i=1}^{N_w} 
\kappa^2 \langle v_{i,k} q_{ik}, \tilde{m} \rangle .
\end{align*}
This system of equations should be compared to the system of equations for the
Hessian~\eqref{eq:Hz}.  From this, it should be clear that the computation
of~$m^*$ corresponds to the solution of another Hessian system with a right-hand
side depending on the state and adjoint variables,~$\{u_i\}$ and~$\{p_i\}$, the
incremental state and adjoint variables,~$\{v_{i,k}\}$ and~$\{q_{i,k}\}$, the
medium parameter~$m$, and the $\{y_k\}$. We denote this right-hand side as
$\mathscr{F}$. In strong form, $m^*$ thus solves
\[ \Hess (m_\map) m^* =
\mathscr{F}(\{u_i\},\{p_i\},\{v_{i,k}\},\{q_{i,k}\},m,\{y_k\}) . \]

\end{document}